%% file: slag_080705.tex
\input amstex
\let\myfrac=\frac
\input eplain.tex
\let\frac=\myfrac
\input epsf




\loadeufm \loadmsam \loadmsbm
\message{symbol names}\UseAMSsymbols\message{,}

\font\myfontdefault=cmr10

\font\mytdmchapfont=cmb10 at 14pt
\font\mytdmheadfont=cmb10 at 10pt
\font\mytdmsubheadfont=cmr10

\magnification 1200
\newif\ifinappendices
\newif\ifundefinedreferences
\newif\ifchangedreferences
\newif\ifloadreferences
\newif\ifmakebiblio
\newif\ifmaketdm

\undefinedreferencesfalse
\changedreferencesfalse


\loadreferencestrue
\makebibliofalse
\maketdmfalse

\def\headpenalty{-400}     
\def\proclaimpenalty{-200} 

%
%

\def\alphanum#1{\ifcase #1 _\or A\or B\or C\or D\or E\or F\or G\or H\or I\or J\or K\or L\or M\or N\or O\or P\or Q\or R\or S\or T\or U\or V\or W\or X\or Y\or Z\fi}
\def\gobbleeight#1#2#3#4#5#6#7#8{}

\newwrite\references
\newwrite\tdm
\newwrite\biblio

\newcount\chapno
\newcount\headno
\newcount\subheadno
\newcount\procno
\newcount\figno
\newcount\citationno

\def\setcatcodes{%
\catcode`\!=0 \catcode`\\=11}%

\ifloadreferences
    {\catcode`\@=11 \catcode`\_=11%
    \input references.tex %
    }%
\else
    \openout\references=references.tex
\fi

\newcount\newchapflag 
\newcount\showpagenumflag 

\global\chapno = -1 
\global\citationno=0
\global\headno = 0
\global\subheadno = 0
\global\procno = 0
\global\figno = 0

\def\resetcounters{%
\global\headno = 0%
\global\subheadno = 0%
\global\procno = 0%
\global\figno = 0%
}

\global\newchapflag=0 
\global\showpagenumflag=0 

\def\chinfo{\ifinappendices\alphanum\chapno\else\the\chapno\fi}%
\def\headinfo{\ifinappendices\alphanum\headno\else\the\headno\fi}
\def\subheadinfo{\the\headno.\the\subheadno}
\def\procinfo{\ifinappendices\alphanum\headno\else\the\headno\fi.\the\procno}
\def\figinfo{\the\headno.\the\figno}
\def\citationinfo{\the\citationno}%
\def\nextheadno{\global\advance\headno by 1 \global\subheadno = 0 \global\procno = 0}
\def\nextsubheadno{\global\advance\subheadno by 1}
\def\nextprocno{\global\advance\procno by 1 \procinfo}
\def\nextfigno{\global\advance\figno by 1 \figinfo}

{\global\let\noe=\noexpand%
%
%
\catcode`\@=11%
\catcode`\_=11%
\setcatcodes%
!global!def!_@@internal@@makeref#1{%
!global!expandafter!def!csname #1ref!endcsname##1{%
!csname _@#1@##1!endcsname%
!expandafter!ifx!csname _@#1@##1!endcsname!relax%
    !write16{#1 ##1 not defined - run saving references}%
    !undefinedreferencestrue%
!fi}}%
!global!def!_@@internal@@makelabel#1{%
!global!expandafter!def!csname #1label!endcsname##1{%
!edef!temptoken{!csname #1info!endcsname}%
!ifloadreferences%
    !expandafter!ifx!csname _@#1@##1!endcsname!relax%
        !write16{#1 ##1 not hitherto defined - rerun saving references}%
        !changedreferencestrue%
    !else%
        !expandafter!ifx!csname _@#1@##1!endcsname!temptoken%
        !else
            !write16{#1 ##1 reference has changed - rerun saving references}%
            !changedreferencestrue%
        !fi%
    !fi%
!else%
    !expandafter!edef!csname _@#1@##1!endcsname{!temptoken}%
    !edef!textoutput{!write!references{\global\def\_@#1@##1{!temptoken}}}%
    !textoutput%
!fi}}%
!global!def!makecounter#1{!_@@internal@@makelabel{#1}!_@@internal@@makeref{#1}}%
!unsetcatcodes%
}
\makecounter{ch}%
\makecounter{head}%
\makecounter{subhead}%
\makecounter{proc}%
\makecounter{fig}%
\makecounter{citation}%
\def\newref#1#2{%
\def\temptext{#2}%
\edef\bibliotextoutput{\expandafter\gobbleeight\meaning\temptext}%
\global\advance\citationno by 1\citationlabel{#1}%
\ifmakebiblio%
    \edef\fileoutput{\write\biblio{\noindent\hbox to 0pt{\hss$[\the\citationno]$}\hskip 0.2em\bibliotextoutput\medskip}}%
    \fileoutput%
\fi}%
\def\cite#1{%
$[\citationref{#1}]$%
\ifmakebiblio%
    \edef\fileoutput{\write\biblio{#1}}%
    \fileoutput%
\fi%
}%
%
%
%

\let\mypar=\par


\def\raggedleft{\leftskip=0pt plus 1fil \parfillskip=0pt}


\font\lettrinefont=cmr10 at 28pt
\def\lettrine #1[#2][#3]#4%
{\hangafter -#1 \hangindent #2
\noindent\hskip -#2 \vtop to 0pt{
\kern #3 \hbox to #2 {\lettrinefont #4\hss}\vss}}

\font\mylettrinefont=cmr10 at 28pt
\def\mylettrine #1[#2][#3][#4]#5%
{\hangafter -#1 \hangindent #2
\noindent\hskip -#2 \vtop to 0pt{
\kern #3 \hbox to #2 {\mylettrinefont #5\hss}\vss}}


\edef\Pagetitle={Blank}

\headline={\hfil\Pagetitle\hfil}

\footline={{\myfontdefault \hfil\folio\hfil}}

\def\nextoddpage
{
\newpage%
\ifodd\pageno%
\else%
    \global\showpagenumflag = 0%
    \null%
    \vfil%
    \eject%
    \global\showpagenumflag = 1%
\fi%
}


\def\newchap#1#2%
{%
%
%
\global\advance\chapno by 1%
\resetcounters%
%
%
\newpage%
\ifodd\pageno%
\else%
    \global\showpagenumflag = 0%
    \null%
    \vfil%
    \eject%
    \global\showpagenumflag = 1%
\fi%
\global\newchapflag = 1%
\global\showpagenumflag = 1%
%
%
{\font\chapfontA=cmsl10 at 30pt%
\font\chapfontB=cmsl10 at 25pt%
\null\vskip 5cm%
{\chapfontA\raggedleft\hfil%
{%
\ifnum\chapno=0
    \phantom{%
    \ifinappendices%
        Annexe \alphanum\chapno%
    \else%
        \the\chapno%
    \fi}%
\else%
    \ifinappendices%
        Annexe \alphanum\chapno%
    \else%
        \the\chapno%
    \fi%
\fi%
}%
\par}%
\vskip 2cm%
{\chapfontB\raggedleft%
\lineskiplimit=0pt%
\lineskip=0.8ex%
\hfil #1\par}%
\vskip 2cm%
}%
\edef\Pagetitle{#2}%
%
%
\ifmaketdm%
    \def\temp{#2}%
    \def\tempbis{\nobreak}%
    \edef\chaptitle{\expandafter\gobbleeight\meaning\temp}%
    \edef\mynobreak{\expandafter\gobbleeight\meaning\tempbis}%
    \edef\textoutput{\write\tdm{\bigskip{\noexpand\mytdmchapfont\noindent\chinfo\ - \chaptitle\hfill\noexpand\folio}\par\mynobreak}}%
\fi%
\textoutput%
}


\def\newhead#1%
{%
\ifhmode%
    \mypar%
\fi%
\ifnum\headno=0%
    \ifinappendices
         \nobreak\vskip -\lastskip%
         \penalty\headpenalty\vskip .5cm%
    \else
    \fi
\else%
    \nobreak\vskip -\lastskip%
    \penalty\headpenalty\vskip .5cm%
\fi%
\nextheadno%
\ifmaketdm%
    \def\temp{#1}%
    \edef\sectiontitle{\expandafter\gobbleeight\meaning\temp}%
    \edef\textoutput{\write\tdm{\noindent{\noexpand\mytdmheadfont\quad\headinfo\ - \sectiontitle\hfill\noexpand\folio}\par}}%
    \textoutput%
\fi%
\font\headfontA=cmbx10 at 14pt%
{\headfontA\noindent\headinfo\ -\ #1.\hfil}%
\nobreak\vskip .5cm%
}%


\def\newsubhead#1%
{%
\ifhmode%
    \mypar%
\fi%
\ifnum\subheadno=0%
\else%
    \penalty\headpenalty\vskip .4cm%
\fi%
\nextsubheadno%
\ifmaketdm%
    \def\temp{#1}%
    \edef\subsectiontitle{\expandafter\gobbleeight\meaning\temp}%
    \edef\textoutput{\write\tdm{\noindent{\noexpand\mytdmsubheadfont\quad\quad\subheadinfo\ - \subsectiontitle\hfill\noexpand\folio}\par}}%
    \textoutput%
\fi%
\font\subheadfontA=cmsl10 at 12pt
{\subheadfontA\noindent\subheadinfo\ #1.\hfil}%
\nobreak\vskip .25cm%
}%

%
%


\font\mathromanten=cmr10
\font\mathromanseven=cmr7
\font\mathromanfive=cmr5
\newfam\mathromanfam
\textfont\mathromanfam=\mathromanten
\scriptfont\mathromanfam=\mathromanseven
\scriptscriptfont\mathromanfam=\mathromanfive
\def\mathroman{\fam\mathromanfam}


\font\sansseriften=cmss10
\font\sansserifseven=cmss7
\font\sansseriffive=cmss5
\newfam\sansseriffam
\textfont\sansseriffam=\sansseriften
\scriptfont\sansseriffam=\sansserifseven
\scriptscriptfont\sansseriffam=\sansseriffive
\def\mathsf{\fam\sansseriffam}


\font\boldten=cmb10
\font\boldseven=cmb7
\font\boldfive=cmb5
\newfam\mathboldfam
\textfont\mathboldfam=\boldten
\scriptfont\mathboldfam=\boldseven
\scriptscriptfont\mathboldfam=\boldfive
\def\mathbf{\fam\mathboldfam}


\font\mycmmiten=cmmi10
\font\mycmmiseven=cmmi7
\font\mycmmifive=cmmi5
\newfam\mycmmifam
\textfont\mycmmifam=\mycmmiten
\scriptfont\mycmmifam=\mycmmiseven
\scriptscriptfont\mycmmifam=\mycmmifive

\def\hexa#1{\ifcase #1 0\or 1\or 2\or 3\or 4\or 5\or 6\or 7\or 8\or 9\or A\or B\or C\or D\or E\or F\fi}
\mathchardef\mathi="7\hexa\mycmmifam7B
\mathchardef\mathj="7\hexa\mycmmifam7C


\font\mymsbmten=msbm10 at 8pt
\font\mymsbmseven=msbm7 at 5.6pt
\font\mymsbmfive=msbm5 at 4pt
\newfam\mymsbmfam
\textfont\mymsbmfam=\mymsbmten
\scriptfont\mymsbmfam=\mymsbmseven
\scriptscriptfont\mymsbmfam=\mymsbmfive

\mathchardef\mybeth="7\hexa\mymsbmfam69
\mathchardef\mygimmel="7\hexa\mymsbmfam6A
\mathchardef\mydaleth="7\hexa\mymsbmfam6B


\def\placelabel[#1][#2]#3{{%
\setbox10=\hbox{\raise #2cm \hbox{\hskip #1cm #3}}%
\ht10=0pt%
\dp10=0pt%
\wd10=0pt%
\box10}}%


\newif\ifinproclaim%
\global\inproclaimfalse%
\def\proclaim#1{%
\medskip%
%
%
\bgroup%
\inproclaimtrue%
\setbox10=\vbox\bgroup\leftskip=0.8em\noindent{\bf #1}\sl%
}

\def\endproclaim{%
\egroup%
\setbox11=\vtop{\noindent\vrule height \ht10 depth \dp10 width 0.1em}%
\wd11=0pt%
\setbox12=\hbox{\copy11\kern 0.3em\copy11\kern 0.3em}%
\wd12=0pt%
\setbox13=\hbox{\noindent\box12\box10}%
\noindent\unhbox13%
\egroup%
\medskip\ignorespaces%
}

\def\proclaim#1{%
\medskip%
\bgroup%
\inproclaimtrue%
\noindent{\bf #1}%
\nobreak\medskip%
\sl%
}

\def\endproclaim{%
\mypar\egroup\penalty\proclaimpenalty\medskip\ignorespaces%
}

\def\noskipproclaim#1{%
\medskip%
\bgroup%
\inproclaimtrue%
\noindent{\bf #1}\nobreak\sl%
}

\def\endnoskipproclaim{%
\mypar\egroup\penalty\proclaimpenalty\medskip\ignorespaces%
}


\def\minn{{m\in\Bbb{N}}}
\def\ninn{{n\in\Bbb{N}}}

\def\proof{{\noindent\bf Proof:\ }}

\def\remark{{\noindent\sl Remark:\ }}
\def\suite#1#2{({#1}_{#2})_{#2\in\Bbb{N}}}

\def\msf#1{{\mathsf #1}}

\def\qed{~$\square$}
\def\munion{\mathop{\cup}}
\def\minter{\mathop{\cap}}
\def\myitem#1{%
\ifinproclaim%
    \item{#1}%
\else%
    \noindent\hbox to .5cm{\hfill#1\hss}
\fi}%

\catcode`\@=11
\def\Eqalign#1{\null\,\vcenter{\openup\jot\m@th\ialign{%
\strut\hfil$\displaystyle{##}$&$\displaystyle{{}##}$\hfil%
&&\quad\strut\hfil$\displaystyle{##}$&$\displaystyle{{}##}$%
\hfil\crcr #1\crcr}}\,}
\catcode`\@=12

\def\makeop#1{%
\global\expandafter\def\csname op#1\endcsname{{\mathroman #1}}}%

\def\makeopsmall#1{%
\global\expandafter\def\csname op#1\endcsname{{\mathroman{\lowercase{#1}}}}}%

\makeopsmall{ArcTan}%
\makeopsmall{ArcCos}%
\makeop{Arg}%
\makeop{Det}%
\makeop{Log}%
\makeop{Re}%
\makeop{Im}%
\makeop{Dim}%
\makeopsmall{Tan}%
\makeop{Ker}%
\makeopsmall{Cos}%
\makeopsmall{Sin}%
\makeop{Exp}%
\makeopsmall{Tanh}%
\makeop{Tr}%
\makeop{End}%
\makeop{Long}%
\makeop{Ch}%
\makeop{Exp}%
\makeop{Int}%
\makeop{Ext}%
\makeop{Aire}%
\makeop{Im}%
\makeop{Conf}%
\makeop{Exp}%
\makeop{Mod}%
\makeop{Log}%
\makeop{Ext}%
\makeop{Int}%
\makeop{Dist}%
\makeop{Aut}%
\makeop{Id}%
\makeop{SO}%
\makeop{Homeo}%
\makeop{Vol}%
\makeop{Ric}%
\makeop{Hess}%
\makeop{Euc}%
\makeop{Isom}%
\makeop{Max}%
\makeop{Long}%
\makeop{Fixe}%
\makeop{Wind}%
\makeop{Mush}%
\makeop{Ad}%
\makeop{loc}%
\makeop{Len}%
\makeop{Area}%
\makeop{SL}%
\makeop{GL}%
\makeop{dVol}%
\makeop{Min}%
\makeop{Symm}%
\makeop{O}%
\makeop{Cosh}

\let\emph=\bf

\hyphenation{quasi-con-formal}

%
%

\ifmakebiblio%
    \openout\biblio=biblio.tex %
    {%
        \edef\fileoutput{\write\biblio{\bgroup\leftskip=2em}}%
        \fileoutput
    }%
\fi%

\newref{Aronsz}{Aronszajn N., A unique continuation theorem for elliptic differential equations or inequalities of the second order, {\sl J. Math. Pures Appl.} {\bf 36} (1957),
2359--239}
\newref{BallGromSch}{Ballman W., Gromov M., Schroeder V., {\it Manifolds of nonpositive curvature}, Progress in Mathematics, {\bf 61}, Birkh\"auser, Boston, (1985)}
\newref{CalabiB}{Calabi E., Improper affine hyperspheres of convex type and a generalisation of a theorem by K. J\"orgens, {\sl Michigan Math. J.} {\bf 5} (1958), 105--126} 
\newref{Corlette}{Corlette K., Immersions with bounded curvature, {\sl Geom. Dedicata} {\bf 33} (1990), no. 2, 153--161}
\newref{GilbTrud}{Gilbard D., Trudinger N. S., {\sl Elliptic partical differential equations of second order}, Die Grund\-lehren der mathemathischen Wissenschagten, {\bf 224},
Springer-Verlag, Berlin, New York (1977)}
\newref{GromB}{Gromov M., Pseudoholomorphic curves in symplectic manifolds, {\sl Invent. Math.} {\bf 82} (1985),  no. 2, 307--347}
\newref{GromA}{Gromov M., {\sl Metric Structures for Riemannian and Non-Riemannian Spaces}, Progress in Mathematics, {\bf 152}, Birkh\"auser, Boston, (1998)}
\newref{HarveyLawson}{Harvey R., Lawson H. B. Jr., Calibrated geometries, {\sl Acta. Math.} {\bf 148} (1982), 47--157}
\newref{Hitchin}{Hitchin N., The moduli space of special Lagrangian submanifolds, {\sl Ann. Scuola Norm. Sup. Pisa Cl. Sci. (4)} {\bf 25} (1997), no. 3-4, 503--515}
\newref{JostXin}{Jost J., Xin Y.L., A Bernstein theorem for special Lagrangian graphs, {\sl Calc. Var.} {\bf 15} (2002), 299--312}
\newref{Joyce}{Joyce D., Lectures on special Lagrangian geometry, {\sl Global theory of minimal surfaces}, Clay Math. Proc., Amer. Math. Soc., Providence, RI, (2005), {\bf 2}, 667--695.}
\newref{LabC}{Labourie F., Probl\`eme de Minkowksi et surfaces \`a courbure constante dans les vari\'e\t'es hyperboliques, {\sl Bull. Soc. Math. France} {\bf 119} (1991), no. 3, 307--325}
\newref{LabB}{Labourie F., Probl\`emes de Monge-Amp\`ere, courbes holomorphes et laminations, {\sl GAFA} {\bf 7}, no. 3, (1997), 496--534}
\newref{LabA}{Labourie F., Un lemme de Morse pour les surfaces convexes, {\sl Invent. Math.} {\bf 141} (2000), 239--297}
\newref{Peterson}{Peterson P., {\sl Riemannian Geometry}, Graduate Texts in Mathematics, {\bf 171}, Springer Verlag, New York, (1998)}
\newref{Pog}{Pogorelov A. V., On the improper convex affine hyperspheres, {\sl Geometriae Dedicata} {\bf 1} (1972), 33--46}
\newref{RosSpruck}{Rosenberg H., Spruck J. On the existence of convex hyperspheres of constant Gauss curvature in hyperbolic space, {\sl J. Diff. Geom.} {\bf 40} (1994), no. 2,
379--409} 
\newref{SalMcDuffB}{McDuff D., Salamon D., {\sl $J$-holomorphic curves and quantum cohomology}, University Lecture Series, {\bf 6}, AMS, Providence, (1994)}
\newref{SalMcDuffA}{McDuff D., Salamon D., {\sl Introduction to symplectic topology}, Oxford, (1995)}
\newref{SmiE}{Smith G., Th\`ese de doctorat, Paris (2004)}%
\newref{SmiG}{Smith G., Pointed $k$-surfaces, {\sl Bull. Soc. Math. France} {\bf 134} (2006), no. 4, 509--557}
\newref{SmiF}{Smith G., An Arzela-Ascoli theorem for immersed submanifolds, In Preparation}%
\newref{SmiH}{Smith G., Constant special Lagrangian hypersurfaces in quasi-Fuchsian hyperbolic ends, in preparation}
\newref{Smoczyk}{Smoczyk K., Longtime existence of the Lagrangian mean curvature flow, {\sl Calc. Var. Partial Differential Equations} {\bf 20} (2004), no. 1, 25--46}
\newref{Yuan}{Yuan Y., A Berstein problem for special Lagrangian equations, {\sl Invent. Math.} {\bf 150} (2002), 117--125}

\ifmakebiblio%
    {\edef\fileoutput{\write\biblio{\egroup}}%
    \fileoutput}%
\fi%

%
%
%
\document
\myfontdefault
\global\chapno=1
\global\showpagenumflag=1
\def\Pagetitle{}
\null
\vfill
\def\centre{\rightskip=0pt plus 1fil \leftskip=0pt plus 1fil \spaceskip=.3333em \xspaceskip=.5em \parfillskip=0em \parindent=0em}%
\def\textmonth#1{\ifcase#1\or January\or Febuary\or March\or April\or May\or June\or July\or August\or September\or October\or November\or December\fi}
\font\abstracttitlefont=cmr10 at 14pt
{\abstracttitlefont\centre Special Lagrangian Curvature\par}
\bigskip
{\centre Graham Smith\par}
\bigskip
{\centre \the\day\ \textmonth\month\ \the\year\par}
\bigskip
{\centre Max Planck Institute,\par
Vivatsgasse 7.,\par
53111 Bonn,\par
GERMANY\par}
\bigskip
\noindent{\emph Abstract:\ } We define the notion of special Lagrangian curvature, showing how it may be interpreted as an alternative higher
dimensional generalisation of two dimensional Gaussian curvature. We obtain first a local rigidity result for this curvature when the ambiant manifold
has negative sectional curvature. We then show how this curvature relates to the canonical special Legendrian structure of spherical subbundles of the 
tangent bundle of the ambiant manifold. This allows us to establish a strong compactness result. In the case where the special Lagrangian angle equals 
$(n-1)\pi/2$, we obtain compactness modulo a unique mode of degeneration, where a sequence of hypersurfaces wraps ever 
tighter round a geodesic.
\par
\vfill
\nextoddpage
\def\Pagetitle{\sl Special Lagrangian Curvature}
\global\pageno=1
\newhead{Introduction}
\newsubhead{Curvature and Weingarten Problems} 
\noindent Gaussian curvature is a fundamental object of study in differential geometry. Despite being simple to express, it leads to very hard
problems involving highly non-linear PDEs. In \cite{LabB} and \cite{LabA}, by showing how in two dimensions 
surfaces of constant Gaussian curvature may be studied in terms of pseudo-holomorphic curves in the unitary bundle of the ambiant manifold, Labourie 
used the powerful machinery first developed by
Gromov in \cite{GromB} to obtain a number of elegant results concerning the existence and structure of constant Gaussian curvature surfaces in
negatively curved spaces. Amongst other things, he applies these results in \cite{LabC} to the study of geometrically finite, three dimensional, hyperbolic manifolds.
These ideas were also used by the author in \cite{SmiG} to obtain a satisfying result relating geometrically finite surfaces of constant Gaussian curvature in
hyperbolic space to holomorphic ramified coverings of the Riemann sphere.
\medskip
\noindent In higher dimensions, the difficulty is aggravated by the loss of uniform ellipticity of the PDE defining Gaussian curvature. Nonetheless, by leading a frontal
attack on this PDE, Rosenberg and Spruck were able to obtain in \cite{RosSpruck} strong existence results for complete hypersurfaces of constant Gaussian curvature immersed 
in hyperbolic space and satisfying certain boundary conditions. What the compactness properties of families of such hypersurfaces may be remains, however, an open question.
\medskip
\noindent This paper arose from the attempt to generalise the techniques developed by Labourie to higher dimensions. We found that an alternative higher dimensional generalisation
of Gaussian curvature yields a much more tractable problem. We thus obtain what we have chosen to call special Lagrangian curvature. This curvature is intimately related
to the canonical special Legendrian structure of the unitary bundle (see the following section), and reduces in special cases to elementary combinations of more classical
curvatures (see below).
\medskip
\noindent This curvature may be defined in various ways. In this introduction, we present it in a form which is only defined for convex hypersurfaces. Although
this definition is slightly technical, it exhibits the most clearly its geometric significance. Let $A$ be a positive definite, symmetric matrix. Let $\lambda_1,...,\lambda_n$ be
the eigenvalues of $A$. For $r>0$, we define $SL_r(A)$ by:
$$
SL_r(A) = \opArg(\opDet(\opId + i r A)) = \sum_{i=1}^n \opArcTan(r \lambda_i).
$$
\noindent $SL_r$ is a strictly increasing function of $r$. Moreover, $SL_0=0$ and $SL_\infty=n\pi/2$. Thus, for all $\theta\in]0,n\pi/2[$, there exists a unique $r>0$ such that:
$$
SL_r(A) = \theta.
$$
\noindent We now define $R_\theta(A) = r$. $R_\theta$ is invariant under the action of $O(n)$ on the space of positive definite, symmetric matrices, and it may thus be used
to define curvature.
\medskip
\noindent Let $(M,g)$ be a Riemannian manifold of dimension $(n+1)$, and let $\Sigma=(S,i)$ be an oriented immersed hypersurface in $M$ (ie. an immersed submanifold of codimension $1$). 
let $A$ be the shape (Weingarten) operator of $\Sigma$. We say that $\Sigma$ is convex if and only if $A$ is positive definite. When $\Sigma$ is convex, for $\theta\in ]0,n\pi/2[$, we 
define $R_\theta(\Sigma)$, the $\theta$-special Lagrangian curvature of $\Sigma$ by:
$$
R_\theta(\Sigma) = R_\theta(A).
$$
\noindent Of most interest is the case when $\theta$ is in the half-open interval $[(n-1)\pi/2,n\pi/2[$. Here, {\it convexity is naturally related to the curvature} (for lower values
of $\theta$, a sequence of hypersurfaces of constant special Lagrangian curvature can degenerate by ceasing to be convex). In particular, the case $\theta=(n-1)\pi/2$ yields
the most interesting geometry. Here the special Lagrangian curvature has the simplest form, and here the
compactness result yields a fascinating form of degeneration where 
hypersurfaces wrap ever more tightly round geodesics, thus generalising to higher dimensions the so-called ``curtain surfaces'', studied by Labourie in \cite{LabB}. Examining
the lower dimensional cases illustrates the natural geometric significance of this curvature:
\medskip
\myitem{(i)} when $n=2$, $R_{\pi/2}^{-2} = K$, where $K$ is the Gaussian curvature of $\Sigma$;
\medskip
\myitem{(ii)} when $n=3$, $R_{\pi}^{-2} = K/H$, where $H$ is the mean curvature of $\Sigma$; and
\medskip
\myitem{(iii)} in general $R_{(n-1)\pi/2}$ is the unique value of $r$ such that:
$$
\chi_n - r^2\chi_{n-2} + r^4\chi_{n-4} - ... = 0,
$$
\noindent where, for all $i$, we define the $i$'th higher principal curvature $\chi_i$ of $\Sigma$ such that, for all $t\in\Bbb{R}$:
$$
\opDet(I+tA) = \sum_{i=0}^n\chi_i(A)t^i.
$$
\noindent We thus see that the study of special Lagrangian curvature neatly forms a special case of the study of Weingarten hypersurfaces.
\medskip
\noindent Although the form $R_\theta$ of the special Lagrangian curvature is more transparent to geometry, the form $SL_r$ is much more tractable to analysis. Trivially, $SL_r$ is
constant and equal to $\theta$ if and only if $R_\theta$ is constant and equal to $r$. Therefore, in the sequel, we work with the $r$-special Lagrangian curvature, defined to be
equal to $SL_r(A)$, and we will study hypersurfaces of constant $r$-special Lagrangian curvature.
\medskip
\noindent The main results of this paper provide two key tools for the study of these hypersurfaces. The first is {\it local rigidity} (Theorem \procref{ThmRigidity})
in the case where the ambiant manifold has negative sectional curvature, and the second, and by far the more significant, is {\it precompactness} (Theorem 
\procref{ThmCompactnessOfHypersurfaces}).
\medskip
\noindent As an illustration of the precompactness result, in the very special case where $n=3$ and $\theta=\pi$, a direct application of Theorem \procref{ThmCompactnessOfHypersurfaces}
yields:
\goodbreak
\proclaim{Theorem \nextprocno}
\noindent Let $M$ be a Riemanian manifold of dimension $4$. Let $\kappa>0$ be a positive real number. Let $(\Sigma_n,p_n)_\ninn$ be a family of pointed, immersed hypersurfaces in $M$ and let
$K_n$ and $H_n$ be the Gaussian and mean curvatures respectively of $\Sigma_n$. Suppose that, for all $n$, $K_n/H_n=\kappa$, then, either:
\medskip
\myitem{(i)} $(\Sigma_n,p_n)_\ninn$ subconverges smoothly to a pointed, immersed hypersurface $(\Sigma_0,p_0)$ in $M$; or
\medskip
\myitem{(ii)} $(\Sigma_n,p_n)_\ninn$ contains a subsequence which degenerates by converging to a complete geodesic.
\endproclaim
\remark Large families of examples of such hypersurfaces are described in section $3$ as well as the forthcoming paper \cite{SmiH}.
\medskip
\remark The precise modes of convergence in both cases are described explicitely in Theorem \procref{ThmCompactnessOfHypersurfaces}.
\medskip
\remark This result suggests that $K/H$ be considered as an alternative natural generalisation of two dimensional Gaussian
curvature in the sense that it exhibits identical geometric behaviour to that described for Gaussian curvature by Labourie in \cite{LabB}.
\newsubhead{Positive Special Legendrian Structures}
\noindent The results of this paper also apply in a very different setting: that of positive special Legendrian structures. Special Legendrian structures are the contact equivalent
of special Lagrangian structures, which were first introduced in the landmark paper \cite{HarveyLawson} of Harvey and Lawson and have since been of considerable interest to mathematicians
and theoretical physicists (one excellent survey may be found in \cite{Joyce} and a thoroughly enjoyable read in \cite{Hitchin}). We believe that, in this setting, our results will have
wider applications beyond the scope of this paper.
\medskip
\noindent We define a positive special Legendrian structure are follows: let $(M,g)$ be a Riemannian manifold of dimension $(2n+1)$. Let $\Cal{P}$ be a principal $\opSO(n)$ bundle over
$M$. We define the action $\varphi$ of $\opSO(n)$ over $\Bbb{R}^n\oplus\Bbb{R}^n$ by:
$$
\varphi(A)(u,v) = (Au,Av).
$$
\noindent Let $\omega_n$ be the canonical symplectic structure over $\Bbb{R}^n\oplus\Bbb{R}^n$. Let $F$ be the associated bundle over $M$ defined by the representation $\varphi$:
$$
F = \Cal{P}\otimes_{\varphi(\opSO(n))}\Bbb{R}^n\oplus\Bbb{R}^n.
$$
\noindent Since it is preserved by the action of $\opSO(n)$, $\omega_n$ defines a form $\hat{\omega}_n$ over $F$. We define a positive special Legendrian structure over $M$ to be a triple $(\Cal{P},W,\Phi)$ where:
\medskip
\myitem{(i)} $\Cal{P}$ is an $\opSO(n)$ principal bundle over $M$,
\medskip
\myitem{(ii)} $W$ is a contact structure over $M$, and
\medskip
\myitem{(iii)} $\Phi:F\rightarrow W$ is an isomorphism of vector bundles,
\medskip
\noindent such that if we denote by $\omega$ the canonical symplectic form over $W$ then $\omega=\Phi_*\hat{\omega}_n$. In fact, it suffices to assume that 
$\Phi_*\hat{\omega}_n$ and $\omega$ are merely colinear, but we impose the stricter condition of equality for ease of presentation.
\medskip
\noindent We identify $\Bbb{R}^n\oplus\Bbb{R}^n$ with $\Bbb{C}^n$. Let $\Omega_n = dz^1\wedge ... \wedge dz^n$ be the canonical special Lagrangian form over 
$\Bbb{R}^n\oplus\Bbb{R}^n$. Similarly, we define the canonical Minkowski metric $m_n$ over $\Bbb{R}^n\oplus\Bbb{R}^n$ by:
$$
m_n( (X,Y),(X,Y)) = \langle X,Y\rangle.
$$
\noindent Both $\Omega_n$ and $m_n$ are preserved by $\varphi(\opSO(n))$. They consequently induce forms
$\hat{\Omega}_n$ and $\hat{m}_n$ over $F$. Likewise, the canonical metric $g_n$ on $\Bbb{C}^n$ induces a metric $\hat{g}_n$ over $F$. 
Thus, by pushing forward through the isomorphism $\Phi$, we obtain forms over $W$ that we may denote by $\Omega$, $m$ and $g$ respectively. 
The positive special Legendrian structure thus defines a quadruplet of forms $(\Omega,m,\omega,g)$ over the contact bundle $W$. 
Moreover, the metric over $W$ may trivially be extended to a metric over $M$. Alternatively, given such a quadruplet, satisfying 
certain algebraic relations, we may deduce a positive special Legendrian structure over $M$, and we will adopt this latter point of view in the sequel.
\medskip
\noindent The motivation for the introduction of a Minkowski metric arises from the following result of Jost \& Xin and Yuan (\cite{JostXin} and \cite{Yuan}):
\proclaim{Theorem \procref{ThmYuanPlanes}\ {\bf [Jost, Xin, 2002]}, {\bf [Yuan, 2002]}}
\noindent Let $\Sigma$ be a complete immersed special Lagrangian submanifold of $\Bbb{R}^n\oplus\Bbb{R}^n$ of type
$C^{1,\alpha}$. If $\Sigma$ is spacelike with respect to the canonical Minkowski metric $m$ over
$\Bbb{R}^n\oplus\Bbb{R}^n$ (i.e. if the restriction of $m$ to $T\Sigma$ is non-negative), then $\hat{\Sigma}$ is an
affine subspace.
\endproclaim
\remark In fact, the use of the positivity condition in \cite{JostXin} was itself inspired by the work \cite{Smoczyk} of Smoczyk concerning the Lagrangian mean curvature
flow. Smoczyk showed that this condition ensures that the tangent space of the special Lagrangian submanifold remains within a convex subset of the Grassmannian of Lagrangian
subspaces, and it is precisely this property that is used in \cite{JostXin} to prove the theorem.
\medskip
\noindent It is this result, of Bernstein type, which forms the core of our main compactness result, (Theorem \procref{ThmCompactness}), which we reduce to the former by means of
a blow-up type argument.
\newsubhead{Principal Results - Positive Special Legendrian Submanifolds}
\noindent Let $(M,g)$ be a Riemannian manifold of dimension $(2n+1)$ and let $(\Cal{P},W,\Phi)$ be a special Legendrian structure over $M$. We define 
$\Cal{SL}_\theta^+(M)$ to be the family of all complete, pointed, positive special Legendrian submanifolds immersed in $M$ (see section \subheadref{HeadStructuresSLClassiques}). 
That is $(\Sigma,p)=(S,i,p)$ lies in $\Cal{SL}_\theta^+(M)$ if and only if it is complete and:
$$
\omega|_{T\Sigma},\qquad \opIm(e^{-i\theta}\Omega)|_{T\Sigma} = 0,\qquad m|_{T\Sigma}\geqslant 0.
$$
\noindent We define the canonical projection $\pi:\Cal{SL}_\theta^+(M)\rightarrow M$ by:
$$
\pi(\Sigma,p) = p.
$$
\noindent If we provide $\Cal{SL}_\theta^+(M)$ with the Cheeger/Gromov topology (see section 
\subheadref{ImmersedSubmanifoldsAndTheCheegerGromovTopology}) then the projection $\pi$ is automatically continuous. We obtain the following result:
\proclaim{Theorem \nextprocno\ {\bf Precompactness}}
\noindent The canonical projection $\pi:\Cal{SL}_\theta^+\rightarrow M$ is a proper mapping.
\endproclaim
\proclabel{ThmCompactness}
\remark The same techniques may be used to yield an analogous result for positive special Lagrangian submanifolds.
\newsubhead{Principal Results - Constant Special Lagrangian Curvature}
\noindent We now return our attention to the problem of hypersurfaces of constant $r$-special Lagrangian curvature.
\proclaim{Theorem \nextprocno\ {\bf Local Rigidity}}
\noindent Let $M$ be a Riemannian manifold of dimension $(n+1)$ and of negative sectional curvature bounded above by $-1$. Let $\Sigma = (S,i)$ and $\Sigma' = (S,i')$ be 
two compact, convex, immersed hypersurfaces of constant $r$-special Lagrangian curvature equal to $\theta\leqslant n\opArcTan(r)$. If $i'$ is sufficiently close to $i$ in the $C^0$ topology, then $i'=i$.
\endproclaim
\proclabel{ThmRigidity}
\noindent This result may also be used in the construction of hypersurfaces of constant $r$-special Lagrangian curvature by smoothly deforming the metric of the ambiant
manifold. Indeed, rigidity ensures through the Fredholm alternative that, for small deformations, a smooth family of immersions may be constructed that appropriately 
follows these deformations. This method of construction will be used, in particular, in the forthcoming paper \cite{SmiH} to construct immersed hypersurfaces inside hyperbolic
ends.
\medskip
\noindent In order to state the precompactness result for hypersurfaces of constant special Lagrangian curvature, we require some notation. Let $(M,g)$ be a 
Riemannian manifold of dimension $(n+1)$. For $r>0$, we define $S_r M$, the $r$-sphere bundle over $M$ to be the bundle 
of spheres of radius $r$ in $TM$. Let $\Sigma=(S,i)$ be an oriented immersed hypersurface in $M$. Let $\msf{N}$ be the normal vector field over $\Sigma$. We define
$\hat{\Sigma}_r=(S,\hat{\mathi}_r)$, the $r$-Gauss lifting of $\Sigma$ in $S_r M$ by:
$$
\hat{\mathi}_r = r\msf{N}.
$$
\noindent We define $\Cal{F}_{r,\theta}(M)$ to be the family of all pointed submanifolds $(\Sigma,p)$ immersed in $M$, of constant $r$-special Lagrangian curvature equal 
to $\theta$ whose $r$-Gauss lifting on $S_r M$ is complete. We define $\hat{\Cal{F}}_{r,\theta}(M)$ to be the family of all $r$-Gauss liftings of pointed hypersurfaces 
in $\Cal{F}_{r,\theta}(M)$. We provide $\hat{\Cal{F}}_{r,\theta}(M)$ with the pointed Cheeger/Gromov topology.
\medskip
\noindent We show that $S_r M$ may be furnished with a canonical positive special Legendrian structure. We then
show that $\hat{\Cal{F}}_{\rho,\theta}(M)$ is a subset of $\Cal{SL}^+_\theta(S_\rho M)$. Consequently, if we
define the canonical projection $\pi$ over $\hat{\Cal{F}}_{\rho,\theta}(M)$ as before, we may use Theorem
\procref{ThmCompactness} to show that for every compact subset $K$ of $M$, the subset $\pi^{-1}(K)$ of
$\hat{\Cal{F}}_{\rho,\theta}(M)$ is relatively compact in the space of complete, immersed submanifolds in
$S_\rho M$. We obtain the following result concerning the topological boundary of $\hat{\Cal{F}}_{\rho,\theta}(M)$ in $\Cal{SL}^+_\theta(S_\rho M)$:
\goodbreak
\proclaim{Theorem \nextprocno\ {\bf Precompactness}}
\noindent Let $M$ be a Riemannian manifold of dimension $(n+1)$. Let $\rho\in\Bbb{R}^+$ be a positive real number.
\medskip
\myitem{(i)} If $\theta\in ](n-1)\pi/2,n\pi/2[$, then $\partial\hat{\Cal{F}}_{\rho,\theta}(M)$ is empty.
\medskip
\myitem{(ii)} If $\theta=(n-1)\pi/2$ and if $(\Sigma,p)$ is a pointed immersed submanifold in
$\partial\hat{\Cal{F}}_{\rho,(n-1)\pi/2}(M)$, then $\Sigma$ is a normal $\rho$-sphere bundle over a complete geodesic in
$M$.
\endproclaim
\proclabel{ThmCompactnessOfHypersurfaces}
\remark Heuristically, $(ii)$ says that if $(\Sigma_n,p_n)_\ninn$ is a sequence of immersed hypersurfaces of constant $r$-special Lagrangian curvature equal to
$(n-1)\pi/2$, then, either this sequence subconverges to another such hypersurface, or it contains a subsequence that wraps ever more tightly around a complete
geodesic, and this is the only mode of degeneration that may take place. This theorem thus generalises the result \cite{LabB} of Labourie to higher dimensions.
\medskip
\remark The first case, where $\theta\in ](n-1)\pi/2,n\pi/2[$, follows trivially from the fact that if $SL_r(A)$ is not an integer multiple of 
$\pi/2$, then we automatically obtain a-priori upper and lower bounds on $A$. Thus, the shape operator of a hypersurface of constant $r$-special Lagrangian curvature
equal to $\theta$ is bounded, the slope of the Gauss lifting is therefore also bounded, and this ensures that the Gauss lifting can never become vertical in the 
limit.
\medskip
\remark By using the results \cite{Pog} and \cite{CalabiB} of Pogorelov and Calabi instead of those of Jost \& Xin and Yuan, the techniques used to prove 
Theorems \procref{ThmCompactness} and \procref{ThmCompactnessOfHypersurfaces} may be trivially adapted to yield an analogous, though slightly weaker, compactness result 
concerning hypersurfaces of constant Gaussian curvature (see Proposition \procref{CompactnessForGaussianCurvature}).
\newsubhead{Structure of This Paper}
\noindent This paper is structured as follows:
\medskip
\myitem{(i)} In part $2$, we define the concepts used in the sequel, including special Lagrangian curvature and positive special Legendrian structures. We also prove a
useful unique continuation principle.
\medskip
\myitem{(ii)} In part $3$, we prove a weaker version of local rigidity, and we describe how this may be used to construct numerous examples of hypersurfaces of constant 
special Lagrangian curvature.
\medskip
\myitem{(iii)} In part $4$, we use rescaling techniques to prove the compactness result for special Legendrian submanifolds (Theorem \procref{ThmCompactness}).
\medskip
\myitem{(iv)} In part $5$, we prove the second compactness theorem (Theorem \procref{ThmCompactnessOfHypersurfaces}). This is done by first observing that the result can be expressed in terms of the
vanishing of a certain positive function, and then showing that this function is superharmonic, and thus satisfies the maximum principle. This requires rather involved calculations
of laplacians of functions defined over these hypersurfaces. We then use this result to complete the proof of local rigidity (Theorem \procref{ThmRigidity}).
\medskip
\myitem{(v)} In the appendix, we briefly prove a version of the maximum principle required in the proof of Theorem \procref{ThmCompactnessOfHypersurfaces}.
\medskip
\noindent Finally, I would like to thank Mark Haskins, J\"urgen Jost, Fran\c{c}ois Labourie, John Loftin, Pierre Pansu and Jean-Marc Schlenker for their various contributions 
towards the completion of this paper. I would also like to thank the Max Planck Institute for Mathematics in the Sciences in Leipzig and the Max Planck Institute for Mathematics
in Bonn for providing the excellent conditions for the development of this paper in its different stages.
\newhead{Definitions and Unique Continuation}
\newsubhead{Immersed Submanifolds and the Cheeger/Gromov Topology}
\noindent Let $M$ be a smooth Riemannian manifold. An {\emph immersed submanifold\/} is a pair 
$\Sigma=(S,i)$ where $S$ is a smooth manifold and $i:S\rightarrow M$ is a smooth immersion. A 
{\emph pointed immersed submanifold\/} in  $M$ is a pair $(\Sigma,p)$ where $\Sigma=(S,i)$ is an immersed submanifold in $M$ and $p$ is a point in $S$. An {\emph immersed hypersurface\/} is an immersed submanifold
of codimension $1$. We give $S$ the unique Riemannian metric $i^*g$ which makes $i$ into an isometry. We say that $\Sigma$ is {\emph complete\/} if and only if the Riemannian manifold $(S,i^*g)$ is.
\subheadlabel{ImmersedSubmanifoldsAndTheCheegerGromovTopology}
\medskip
\noindent A {\emph pointed Riemannian manifold\/} is a pair $(M,p)$ where $M$ is a Riemannnian manifold and $p$ is a point in $M$. Let $(M_n,p_n)_{n\in\Bbb{N}}$ be a sequence of complete pointed Riemannian manifolds. For all $n$, we denote by $g_n$ the Riemannian metric over $M_n$. We say that the
sequence $(M_n,p_n)_{n\in\Bbb{N}}$ {\emph converges\/} to the complete pointed manifold $(M_0,p_0)$ in the {\emph Cheeger/Gromov topology\/} if and only if for all $n$, there exists a mapping 
$\varphi_n:(M_0,p_0)\rightarrow (M_n,p_n)$, such that, for every compact subset $K$ of $M_0$, there 
exists $N\in\Bbb{N}$ such that for all $n\geqslant N$:
\medskip
\myitem{(i)} the restriction of $\varphi_n$ to $K$ is a $C^\infty$ diffeomorphism onto its image, and
\medskip
\myitem{(ii)} if we denote by $g_0$ the Riemannian metric over $M_0$, then the sequence of metrics $(\varphi_n^*g_n)_{n\geqslant N}$ converges to $g_0$ in the
$C^\infty$ topology over $K$.
\medskip
\noindent We refer to the sequence $(\varphi_n)_{n\in\Bbb{N}}$ as a sequence of {\emph convergence mappings\/} of the sequence $(M_n, p_n)_{n\in\Bbb{N}}$ with respect to the limit $(M_0,p_0)$. The convergence
mappings are trivially not unique. 
\medskip
\noindent Let $(\Sigma_n,p_n)_{n\in\Bbb{N}}=(S_n,p_n,i_n)_{n\in\Bbb{N}}$ be a sequence of complete pointed
immersed submanifolds in $M$. We say that $(\Sigma_n,p_n)_{n\in\Bbb{N}}$ {\emph converges\/} to 
$(\Sigma_0,p_0)=(S_0,p_0,i_0)$ in the {\emph Cheeger/Gromov topology\/} if and only if 
$(S_n,p_n)_{n\in\Bbb{N}}$ converges to $(S_0,p_0)$ in the Cheeger/Gromov topology, and, for every sequence
$(\varphi_n)_{n\in\Bbb{N}}$ of convergence mappings of $(S_n,p_n)_\ninn$ with respect to this limit, and
for every compact subset $K$ of $S_0$, the sequence of functions $(i_n\circ\varphi_n)_{n\geqslant N}$ converges to the function $(i_0\circ\varphi_0)$ in the $C^\infty$ topology
over $K$.
\medskip
\noindent In an analogous manner, for all $k\geqslant 1$ and for all $\alpha$, we may also define the
$C^{k,\alpha}$ Cheeger/Gromov topology for manifolds and immersed submanifolds. In this case, the
convergence mappings are of type $C^{k,\alpha}$ and the metrics converge in the $C^{k-1,\alpha}$ topology over each compact set.
\newsubhead{Special Lagrangian and Legendrian Structures}
\noindent We identify $\Bbb{R}^{2n}$ and $\Bbb{C}^n$. Using the complex coordinate functions, we define various geometric structures. We define the canonical {\emph symplectic form} by:
\headlabel{HeadStructuresSLClassiques}
\subheadlabel{HeadStructuresSLClassiques}
$$
\omega = \frac{i}{2}\sum_{k=1}^n dz^k\wedge d\overline{z}^k.
$$
\noindent We define the canonical {\emph Minkowski metric} by:
$$
m = \opRe\left(\frac{i}{2}\sum_{k=1}^n dz^k\otimes dz^k\right).
$$
\noindent For all $\theta\in\Bbb{R}$, we define the {\emph $\theta$-special Lagrangian} form $\Omega_\theta$ over
$\Bbb{R}^n\oplus\Bbb{R}^n$ by:
$$
\Omega_\theta = e^{-i\theta}dz^1\wedge ... \wedge dz^n.
$$
\noindent The triple $(\omega,m,\Omega_\theta)$ defines a {\emph positive special Lagrangian structure} over
$\Bbb{R}^n\oplus\Bbb{R}^n$. The stabiliser group of this triple in $\opEnd(\Bbb{R}^{2n})$ is the set of all matrices $M$ of the form:
$$
M = \pmatrix B & \cr & B\cr\endpmatrix,
$$
\noindent where $B\in\opSO(n)$. Let $P$ be an $n$ dimensional subspace of $\Bbb{R}^n\oplus\Bbb{R}^n$. $P$ 
is said to be {\emph $\theta$-special Lagrangian\/} if and only if the restrictions of $\omega$ and 
$\Omega_\theta$ to $P$ vanish. $P$ is said to be {\emph positive\/} (or {\emph spacelike\/}) if and only 
if the restriction of $m$ to $P$ is non-negative definite. In the sequel we will write $SL_\theta$ for $\theta$-special Lagrangian and $SL_\theta^+$ for positive $\theta$-special Lagrangian.
\medskip
\noindent Let $M$ be a $(2n+1)$ dimensional manifold bearing a positive special Lagrangian structure as defined in the introduction. Let $\Sigma=(S,i)$ be an immersed submanifold in $M$. We say that $\Sigma$ is {\emph positive $\theta$-special Legendrian\/} if and only if, for all $p\in S$:
\medskip
\myitem{(i)} $T_p\Sigma\subseteq W_{i(p)}$, and
\medskip
\myitem{(ii)} $T_p\Sigma$ is positive $\theta$-special Lagrangian in $W_{i(p)}$.
\medskip
\noindent We define the family $\Cal{SL}_\theta^+(M)$ to be the family of all complete $SL^+_\theta$ pointed submanifolds of $M$. We give this family the Cheeger/Gromov topology.
Trivially, $\Cal{SL}_\theta^+$ forms a closed subset of the set of all complete pointed immersed submanifolds in $M$ with respect to the Cheeger/Gromov topology.
We define the canonical projection $\pi:\Cal{SL}_\theta^+\rightarrow M$ by:
$$
\pi(\Sigma,p) = \pi((S,i),p) = i(p).
$$
\noindent This mapping is trivially continuous.
\newsubhead{The Positive Special Legendrian Structure Over $S_\rho M$}
\noindent Let $M$ be a complete, oriented Riemannian manifold of dimension $(n+1)$. Let $\pi:TM\rightarrow M$ be the canonical projection. Let $VTM$ be the vertical subbundle of $TTM$, and let $HTM$ be the horizontal subbundle of the Levi-Civita connection of $M$. Thus:
$$
TTM = HTM\oplus VTM \cong \pi^*TM\oplus\pi^*TM.
$$
\noindent Using this identification of $TTM$ with $\pi^*TM\oplus\pi^*TM$, for $X,Y\in\Gamma(M,TM)$, we
define $\left\{X,Y\right\}\in\Gamma(TM,TTM)$ by:
$$
i_H\oplus i_V\left\{X,Y\right\} = (\pi^*X,\pi^*Y).
$$
\noindent Every vector field over $TTM$ may be expressed locally in terms of a linear combination of such 
vector fields. For $X,Y,v\in T_pM$, we define $\left\{X,Y\right\}_v\in T_vTM$ in an analogous manner.
\medskip
\noindent For $\rho$ a positive real number, let $S_\rho M$ be the set of vectors of length $\rho$ in $TM$.
We call $S_\rho M$ the $\rho$-sphere bundle over $M$. We define the tautological vector fields $q^H$ and 
$q^V$ over $TM$, by:
$$
q^H(v) = \left\{v,0\right\}_v,\qquad q^V(v) = \left\{0,v\right\}_v.
$$
\noindent Let $HS_\rho M$ and $VS_\rho M$ be the respective restrictions of the bundles $HTM$ and $VTM$ 
to $S_\rho M$. We define $\langle q^H\rangle^\perp$ (resp. $\langle q^V\rangle^\perp$) to be the subspace
orthogonal to $q^H$ (resp. $q^V$) in $HS_\rho M$ (resp. $VS_\rho M$). Trivially:
$$
{i_\rho}_*:TS_\rho M\cong HS_\rho M\oplus\langle q^V\rangle^\perp.
$$
\noindent We define the subbundle $WS_\rho M$ of $TS_\rho M$ by:
$$
WS_\rho M = \langle q^H\rangle^\perp\oplus\langle q^V\rangle^\perp.
$$
\noindent $WS_\rho M$ defines a contact structure over $S_\rho M$. Moreover, if we denote by $\omega$ the 
canonical symplectic form over $WS_\rho M$, then:
$$
\omega(\left\{X_1,X_2\right\},\left\{Y_1,Y_2\right\}) = \langle X_2,Y_1\rangle - \langle X_1,Y_2\rangle.
$$
\noindent In fact, using the metric of $M$, we identify $TM$ and $T^*M$. The canonical Liouville form of $T^*M$ then restricts to a contact form over $S_\rho M$ whose kernel
is $W$. Moreover, the restriction to $W$ of the exterior derivative of this form is colinear with $\omega$. We thus obtain a second construction of $\omega$.
\medskip
\noindent The distributions $\langle q^H\rangle^\perp$, $\langle q^V\rangle^\perp\subseteq TTM$ are canonically isomorphic and inherit canonical metrics and orientations from $TM$. This 
defines a positive special Legendrian structure over $W$. We explicitly construct the forms $\omega$, $m$, $\Omega_\theta$ and $g$. We define $J$ over $\langle q^H\rangle^\perp$ to be the canonical isometry sending this space into $\langle q^V\rangle^\perp$. $J$ then uniquely extends to a complex structure over $W$.
\medskip
\noindent Let $v_p$ be a point in $S_\rho M$. Let $E_1, ..., E_n\in T_pM$ be such that $(E_1, ..., E_n, v_p/\rho)$ forms an oriented, orthonormal basis. We define $X_1, ..., X_n\in \langle q^H\rangle_v^\perp$, and
$Y_1, ..., Y_n\in\langle q^V\rangle_v^\perp$ by:
$$
X_i=\left\{E_i,0\right\}_v,\qquad Y_i=\left\{0,E_i\right\}_v=JX_i\qquad\forall i.
$$
\noindent $X_1,...,X_n,Y_1,...,Y_n$ is an orthonormal basis of $W_v$. Let 
$(X^1, ..., X^n, Y^1, ..., Y^n)$ be the dual basis of $W_v^*$. We define:
$$\matrix
\omega \hfill&= \sum_{i=1}^n X^i\wedge Y^i,\hfill\cr
m \hfill&= \frac{1}{2}\sum_{i=1}^n(X^i\otimes Y^i + Y^i\otimes X^i),\hfill\cr
\Omega_\theta \hfill&= e^{-i\theta}(X^1+iY^1)\wedge...\wedge(X^n+iY^n),\hfill\cr
g \hfill&= \sum_{i=1}^n X^i\otimes X^i + Y^i\otimes Y^i.\hfill\cr
\endmatrix$$
\noindent The triplet $(\omega, m, \Omega_\theta)$ does not depend on $(E_i)_{1\leqslant i\leqslant n}$. 
Since the form $\omega$ coincides with the canonical symplectic form over $W$, it follows that
$(W,\omega,m,\Omega_\theta)$ does indeed define a positive special Legendrian structure over $S_\rho M$.
\newsubhead{Normal Vector Fields, Second Fundamental Form, Convexity}
\noindent Let $M$ be an oriented Riemannian manifold. Let $\Sigma=(S,i)$ be an oriented hypersurface immersed in $M$. Let $\msf{N}_\Sigma\in\Gamma(S,i^*TM)$ be the exterior unit normal of $\Sigma$
in $M$. The sign of $\msf{N}_\Sigma$ depends on the orientations of $S$ and $M$. We define the {\emph Weingarten operator\/} $A_\Sigma:TS\rightarrow TS$ and the {\emph second fundamental form} $II_\Sigma$ of 
$\Sigma$ by:
$$\matrix
A_\Sigma(X) \hfill&= \nabla_X\msf{N}_\Sigma,\hfill\cr
II_\Sigma(X,Y) \hfill&= \langle A_\Sigma X, Y\rangle.\hfill\cr
\endmatrix$$
\noindent $\Sigma$ is said to be {\emph convex\/} at $p\in S$ if and only if the symmetric bilinear form 
$II_\Sigma$ is either positive or negative (but not mixed) at $p$. Through a slight abuse of language, we say that $\Sigma$ is {\emph convex\/} if it is convex at every point. In this case, by reversing the orientation of $S$ if necessary, we may assume in the sequel that $II_\Sigma$ is positive.
\medskip
\noindent For $\rho\in(0,\infty)$, we define the {\emph $\rho$-Gauss lifting\/} $\hat{\Sigma}_\rho=(S,\hat{\mathi}_\rho)$ of $\Sigma$, which is an immersed submanifold of $S_\rho M$, by:
$$
\hat{\Sigma}_\rho = (S,\hat{\mathi}_\rho) = (S, \rho\msf{N}_\Sigma).
$$
\noindent The relationship between hypersurfaces and immersed submanifolds is made clear by the following elementary lemma:
\proclaim{Lemma \nextprocno}
\noindent Let $\Sigma$ be an oriented immersed hypersurface in $M$. $\hat{\Sigma}_\rho$, the $\rho$-Gauss lifting of $\Sigma$, is a Legendrian submanifold of $S_\rho M$. Moreover, if $\Sigma$ is convex, then $\hat{\Sigma}_\rho$ is positive (spacelike) with respect to the Minkowski metric $m$.
\medskip
\noindent Conversely, if $\hat{\Sigma}$ is an immersed Legendrian submanifold in $S_\rho M$ such that
$T\hat{\Sigma}\minter VS_\rho M=0$, then there exists a unique oriented immersed hypersurface $\Sigma$ in $M$ such that $\hat{\Sigma}$ is the $\rho$-Gauss lifting of $\Sigma$. Moreover, if $\hat{\Sigma}$ is positive (spacelike) with respect to the Minkowski metric $m$, then $\Sigma$ is convex.
\endproclaim
\proclabel{LemmaGaussLiftingLegendrian}
\newsubhead{Curvature}
\noindent Let us denote by $\opSymm(\Bbb{R}^n)$ the space of symmetric matrices over $\Bbb{R}^n$. We define $\Phi:\opSymm(\Bbb{R}^n)\rightarrow\Bbb{C}^*$ by:
$$
\Phi(A) = \opDet(I+iA).
$$
\noindent Since $\Phi$ never vanishes and $\opSymm(\Bbb{R}^n)$ is simply connected, there exists a unique analytic function $\tilde{\Phi}:\opSymm(\Bbb{R}^n)\rightarrow\Bbb{C}$
such that:
$$
\tilde{\Phi}(I) = 0,\qquad e^{\tilde{\Phi}(A)} = \Phi(A)\qquad\forall A\in\opSymm(\Bbb{R}^n).
$$
\noindent We define the function $\opArcTan:\opSymm(\Bbb{R}^n)\rightarrow(-n\pi/2,n\pi/2)$ by:
$$
\arctan(A) = \opIm(\tilde{\Phi}(A)).
$$
\noindent This function is trivially invariant under the action of $O(\Bbb{R}^n)$. If $\lambda_1, ..., \lambda_n$ are the eigenvalues of $A$, then:
$$
\opArcTan(A) = \sum_{i=1}^n\opArcTan(\lambda_i).
$$
\noindent Let $M$ be an oriented Riemannian manifold of dimension $n+1$. Let $\Sigma=(S,i)$ be an oriented immersed
hypersurface in $M$. For $\rho\in(0,\infty)$, we define $\opSL_\rho(\Sigma)$, the {\emph $\rho$-special Lagrangian
curvature\/} of $\Sigma$ by:
$$
\opSL_\rho(\Sigma) = \opArcTan(\rho A_\Sigma).
$$
\noindent For $\theta\in [(n-1)\pi/2,n\pi/2[$ and $\rho\in (0,\infty)$, we define the
family $\Cal{F}_{\rho,\theta}(M)$ by:
$$
\Cal{F}_{\rho,\theta}(M) = \left\{ (\Sigma,p)\text{ s.t. }
\matrix
\bullet\ (\Sigma,p)\text{ is an immersed, pointed hypersurface in $M$, }\hfill\cr
\bullet\ \opSL_\rho(\Sigma) = \theta\text{, and}\hfill\cr
\bullet\ \hat{\Sigma}_\rho\text{ is a complete immersed submanifold of }S_\rho M.\hfill\cr
\endmatrix\right\}.
$$
\noindent We observe that hypersurfaces in $\Cal{F}_{\rho,\theta}(M)$ are always convex.
We define $\hat{\Cal{F}}_{\rho,\theta}(M)$ to be the family of $\rho$-Gauss liftings of pointed hypersurfaces in 
$\Cal{F}_{\rho,\theta}(M)$. Trivially, there is a canonical bijections between $\Cal{F}_{\rho,\theta}(M)$ and $\hat{\Cal{F}}_{\rho,\theta}(M)$.
We furnish $\hat{\Cal{F}}_{\rho,\theta}(M)$ with the Cheeger/Gromov topology. Let $\pi_\rho:S_\rho M\rightarrow M$ be the canonical projection. We define the canonical projection $\pi:\hat{\Cal{F}}_{\rho,\theta}(M)\rightarrow M$ by:
$$
\pi(\hat{\Sigma},p) = \pi(S_\rho,\hat{\mathi}_\rho,p) = (\pi_\rho\circ\hat{\mathi}_\rho)(p) = i(p).
$$
\noindent By definition of the Cheeger/Gromov topology, this projection is continuous.
\medskip
\noindent Bearing in mind Lemma \procref{LemmaGaussLiftingLegendrian}, we obtain the following result:
\proclaim{Lemma \nextprocno}
\noindent For all $\rho\in(0,\infty)$ and for all $\theta\in [(n-1)\pi/2,n\pi/2[$:
$$
\hat{\Cal{F}}_{\rho,\theta} \subseteq \Cal{SL}^+_\theta(S_\rho M).
$$
\noindent Conversely, if $(\hat{\Sigma},p)\in\Cal{SL}^+_\theta(S_\rho M)$ and if $T\hat{\Sigma}\minter VS_\rho M=0$, then
$(\hat{\Sigma},p)$ is the $\rho$-Gauss lifting of a convex immersed hypersurface in $M$ and thus:
$$
(\hat{\Sigma},p)\in\hat{\Cal{F}}_{\rho,\theta}.
$$
\endproclaim
\proclabel{EquivalenceSLCurvatureAndSLImmersions}
\noindent The study of convex immersed hypersurfaces of constant $\rho$-special Lagrangian curvature is a refinement of a
special case of the more general theory of Weingarten hypersurfaces. Indeed, the family of hypersurfaces $\Sigma$
satisfying $\opSL_\rho(\Sigma)\in\theta+\pi\Bbb{Z}$ coincides with the family of hypersurfaces $\Sigma$
satisfying:
$$
\opSin(\theta)\sum_{k\geqslant 0}(-1)^k\rho^{2k}\chi_{2k}(\Sigma) - \opCos(\theta)\sum_{k\geqslant 0}(-1)^k\rho^{2k+1}\chi_{2k+1}(\Sigma) = 0,
$$
\noindent where, for all $i$, $\chi_i(\Sigma)=\chi_i(A)$ is the $i$'th higher principal curvature of $\Sigma$ as defined in the introduction. The special Lagrangian curvature 
deserves particular attention since, by Lemma \procref{EquivalenceSLCurvatureAndSLImmersions}, it reflects a uniformly elliptic problem arising from calibrated geometries. This uniform ellipticity permits us to obtain compactness results which are not necessarily valid for Weingarten hypersurfaces in general, even in such apparently simple cases as hypersurfaces of constant Gaussian curvature.
\newsubhead{Unique Continuation}
\noindent The special Lagrangian equation linearises to a Laplacian. We thus obtain the following unique continuation principle:
\proclaim{Lemma \nextprocno}
\noindent Let $(\Sigma,p)=(S,i,p)$ and $(\Sigma',p')=(S',i',p')$ be two special Legendrian submanifolds immersed in $M$. If $i(p)=i'(p')=q$ and if $\Sigma$ and $\Sigma'$ have the 
same $\infty$-jet at $q$, then $(\Sigma,p)$ is locally equivalent to $(\Sigma',p')$ at $q$. In other words, there exist neighbourhoods $U$ and $U'$ of $p$ and $p'$ respectively, and 
a (locally unique) diffeomorphism $\varphi:U\rightarrow U'$ such that $\varphi(p) = p'$ and:
$$
i'\circ\varphi = i.
$$
\endproclaim
\proclabel{ThmUniqueContinuation}
\proof Let $q=i(p)=i'(p')$. Let $\omega_0$ be the canonical symplectic form over $\Bbb{R}^n\times\Bbb{R}^n$. Let $\beta$ be a primitive of $\omega$. Let us define the
contact form $\alpha$ over $\Bbb{R}\times\Bbb{R}^n\times\Bbb{R}^n$ by:
$$
\alpha = dt - \beta.
$$
\noindent By Darboux's theorem (see, for example, \cite{SalMcDuffA}), for an appropriate choice of $\beta$, there exists:
\medskip
\myitem{(i)} a positive real number $\epsilon$,
\medskip
\myitem{(ii)} a neighbourhood $U$ of $q$ in $M$, and
\medskip
\myitem{(iii)} a diffeomorphism $\varphi:(B_\epsilon(0),0)\rightarrow (U,p)$,
\medskip
\noindent such that:
\medskip
\myitem{(i)} $\varphi_*\alpha$ is colinear with the contact structure over $U$, and
\medskip
\myitem{(ii)} $T_0\varphi\cdot(\left\{0\right\}\times\left\{0\right\}\times\Bbb{R}^n) = T_p\Sigma = T_{p'}\Sigma'$.
\medskip
\noindent Since $\Sigma$ and $\Sigma'$ are Legendrian, there exists an open set $\Omega$ of $0$ in $\Bbb{R}^n$ and functions $f,f':\Omega\rightarrow\Bbb{R}$ such that $\varphi^*(\Sigma,p)$ and $\varphi^*(\Sigma',p')$ are locally the graphs over $\Omega$ of $(f,df)$
and $(f',df')$ respectively.
\medskip
\noindent Let $u$ be a function over $\Omega$. There exists a function $M:\Bbb{R}\times\Bbb{R}^n\times\Bbb{R}^n\rightarrow GL(\Bbb{R}^n)$ such that the graph of $(u,du)$
is $\theta$-special Legendrian if and only if:
$$
\opIm(e^{-i\theta}\opDet(I+iM^{-1}(x,u,du)D^2uM(x,u,du))) = 0.
$$
\noindent Differentiating, we obtain, for all $k$, smooth functions $a_k^{ij}$, $b_k^i$ and $c_k$ over $\Bbb{R}^{(n+1)^2}$
such that:
$$\matrix
a_k^{ij}(x,u,Du,D^2u)\partial_i\partial_j(\partial_k u) + b_k^i(x,u,Du,D^2u)\partial_i(\partial_k u)\hfill\cr
\qquad\qquad\qquad + c_k(x,u,Du,D^2u)\partial_k u = 0.\hfill\cr
\endmatrix$$
\noindent Moreover, $a_k^{ij}$ is invertible close to the origin. Since $f$ and $f'$ both satisfy this relation, we have, for all $k$:
$$\matrix
a_k^{ij}(x,f,Df,D^2f)\partial_i\partial_j\partial_k(f-f')\hfill\cr
\qquad + (a_k^{ij}(x,f',Df',D^2f') - a_k^{ij}(x,f,Df,D^2f))\partial_i\partial_j\partial_k f'\hfill\cr
+b_k^{i}(x,f,Df,D^2f)\partial_i\partial_k(f-f')\hfill\cr
\qquad + (b_k^{i}(x,f',Df',D^2f') - b_k^i(x,f,Df,D^2f))\partial_i\partial_k f'\hfill\cr
+c_k(x,f,Df,D^2f)\partial_i\partial_k(f-f')\hfill\cr
\qquad + (c_k(x,f',Df',D^2f') - c_k(x,f,Df,D^2f))\partial_k f'=0.\hfill\cr
\endmatrix$$
\noindent Since $a_k^{ij}$, $b_k^i$ and $c_k$ are smooth, since all the derivatives of $f$ and $f'$ are bounded near $0$, and since $a_k^{ij}$ is uniformly bounded
below near zero, using Taylor's theorem and the preceeding relation, we find that there exists $K\in\Bbb{R}^+$ such that, for all $k$:
$$
\left|\Delta(\partial_k (f-f'))\right| \leqslant K(\|D^2(f-f')\| + \|D(f-f')\| + \|f-f'\|).
$$
\noindent Consider the function $w$ defined by:
$$
w = (f-f', D(f-f')).
$$
\noindent Using the preceeding relation, we find that there exists $L\in\Bbb{R}^+$ such that:
$$
\left|\Delta w\right| \leqslant L(\|Dw\| + \|w\|).
$$
\noindent Since $w$ vanishes to infinite order at the origin, it follows by Aronszajn's unique continuation principle (see \cite{Aronsz} or \cite{SalMcDuffB}) that $w$ vanishes
identically in a neighbourhood of zero. Consequently $f$ coincides with $f'$ in a neighbourhood of zero and the result follows.\qed
\newhead{Local Rigidity and Examples}
\newsubhead{Infinitesimal Deformation of Hypersurfaces}
\noindent Let $M$ be a Riemannian manifold. Let $\opExp:TM\rightarrow M$ be the exponential mapping of $M$. Let $\Sigma=(S,i)$ be an immersed hypersurface of constant
$\rho$-special Lagrangian curvature in $M$. Let $\msf{N}$ be the normal vector field over $\Sigma$ in $M$. For
$f\in C^\infty(\Sigma)$, we define the family of immersed hypersurfaces
$(\Sigma_{f,t})_{t\in(-\epsilon,\epsilon)}=(S,i_{f,t})_{t\in(-\epsilon,\epsilon)}$ by:
$$
i_{f,t}(p) = \opExp_{i(p)}(tf(p)\msf{N}(p)).
$$
\noindent For all $t$, let $\msf{N}_{f,t}$ be the normal vector field over $\Sigma_{f,t}$. For all $t$, let $A_{f,t}$ be the
Weingarten operator of $\Sigma_{f,t}$. In other words:
$$
A_{f,t} u = \nabla_u\msf{N}_{f,t}.
$$
\noindent We define the operator of variation of $\rho$-special Lagrangian curvature $\Cal{L}^\rho_\theta$, such that:
$$
\Cal{L}^\rho_\theta f = \partial_t\opIm(e^{-i\theta}\opDet(I+i\rho A_{f,t}))|_{t=0}.
$$
\noindent We aim to calculate explicitly the operator $\Cal{L}^\rho_\theta$. Let $R$ be the Riemann
curvature tensor of $\nabla$, let $W$ to be the endomorphism of $TS$ defined by:
$$
i_*W(u) = R_{\msf{N_\Sigma} i_*u}\msf{N_\Sigma},
$$
\noindent and let $\opHess(f)$ be the Hessian of the function $f$:
\proclaim{Lemma \nextprocno}
\noindent If $A$ is the Weingarten operator of $\Sigma$, then:
$$\matrix
\Cal{L}^\rho_\theta f \hfill&= \sqrt{I+\rho^2A^2}\bigl(-\opTr((I+\rho^2A^2)^{-1}\rho\opHess(f)) \hfill\cr
&\qquad + f\opTr((I+\rho^2A^2)^{-1}\rho W) - f\opTr((I+\rho^2A^2)^{-1}\rho A^2)\bigr).
\endmatrix$$
\endnoskipproclaim
\proclabel{LemmaInfinitesimalDeformation}
\proof By Proposition $3.1.1$ of \cite{LabA}:
$$
\partial_t A_{f,t}|_{t=0} = fW - \opHess(f) - fA_{f,0}^2.
$$
\noindent Bearing in mind that:
$$
\opIm(e^{-i\theta}\opDet(I+i\rho A_{f,0})) = 0,
$$
\noindent the result now follows by elementary calculation.\qed
\medskip
\noindent It follows that the operator $\Cal{L}^\rho_\theta$ defined over a submanifold of constant 
$\rho$-special Lagrangian curvature equal to $\theta$ is elliptic. Since it acts over the space of real valued functions, it is of index zero, and we would thus like know under which conditions it is injective
(and thus invertible). In \cite{SmiH} we prove:
\proclaim{Lemma \nextprocno}
\noindent Suppose that the sectional curvature of $M$ is less than $-\kappa<0$ and suppose that the $\rho$-special Lagrangian curvature of $\Sigma$ is less than 
or equal to $n\opArcTan(\sqrt{\kappa}\rho)$, then the zeroth order term of the formula given in the preceeding lemma is non-negative:
$$
J=\opTr\left( (I+\rho^2 A^2)^{-1}\rho W\right) - \opTr\left((I+\rho^2 A^2)^{-1}\rho A^2\right) \geqslant 0
$$
\endproclaim
\proclabel{LemmeBijectivite}
\proof This is an elementary, but long and technical exercise involving Lagrange multipliers. See \cite{SmiH} for details.\qed
\medskip
\noindent This yields the following partial version of Theorem \procref{ThmRigidity}
\proclaim{Lemma \nextprocno}
\noindent Let $M$ be a Riemannian manifold of dimension $(n+1)$ and of negative sectional curvature bounded above by $-1$. Let $\Sigma = (S,i)$ and $\Sigma' = (S,i')$ be two compact,
convex immersed hypersurfaces of constant $\rho$-special Lagrangian curvature equal to $\theta\leqslant n\opArcTan(\rho)$. If $i$ is sufficiently close to $i'$ in the
$C^{2,\alpha}$ topology, then (up to reparametrisation) $i'=i$.
\endproclaim
\proclabel{LemmaLocalRigidity}
\remark A complete proof of Theorem \procref{ThmRigidity} will follow immediately from Theorem \procref{ThmCompactnessOfHypersurfaces}.
\medskip
\proof This follows from the preceeding two Lemmata and the implicit function theorem for smooth functions on Banach manifolds.\qed
\newsubhead{Construction of Examples}
\noindent Non trivial examples of hypersurfaces of constant special Lagrangian curvature may be constructed in various ways. Non complete examples may be trivially 
constructed by consider hypersurfaces of revolution and thus transforming the PDE into an ODE.
\medskip
\noindent Another method of construction is the following: let $M$ be a hyperbolic manifold of dimension $(n+1)$. Let $N\subseteq M$ be a totally geodesic hypersurface. For all
$R$, let $N_R$ be the equidistant hypersurface at distance $R$ from $N$. Elementary hyperbolic geometry allows us to calculate:
\noskipproclaim{Lemma \nextprocno}
$$
\opSL_\rho(N_R) = n\opArcTan(\rho\opTanh(R)).
$$
\endnoskipproclaim
\noindent Thus $N_R$ satisfies the hypothesis of Lemma \procref{LemmeBijectivite}. Thus, using Lemmata \procref{LemmaInfinitesimalDeformation} and
\procref{LemmeBijectivite} along with the implicit function theorem for smooth functions on Banach manifolds, we may show that for small, smooth deformations
of the metric on $M$, there exists deformations of $N_R$ whose special Lagrangian curvature remains constant.
\medskip
\noindent In \cite{SmiH}, we use an analogous technique to construct hypersurfaces of constant special Lagrangian curvature in hyperbolic ends.
\newhead{Compactness}
\newsubhead{Preliminary Definitions and Results}
\noindent We collect results required to prove Theorem \procref{ThmCompactness}. First, in \cite{JostXin} and \cite{Yuan}, Jost \& Xin and Yuan prove the following result:
\headlabel{HeadCompacite}
\proclaim{Theorem \nextprocno\ {\bf [Jost, Xin, 2002]}, {\bf [Yuan, 2002]}}
\noindent If $\hat{\Sigma}$ is a complete immersed $SL_\theta^+$ submanifold of $\Bbb{R}^n\oplus\Bbb{R}^n$ of type $C^{1,\alpha}$, then
$\hat{\Sigma}$ is an affine subspace.
\endproclaim
\proclabel{ThmYuanPlanes}
\remark \cite{JostXin} and \cite{Yuan} obtain the same result simultaneously using very different techniques. Jost and Xin use the properties of harmonic maps into convex 
subsets of Grassmannians. Yuan uses geometric measure theory. \cite{JostXin} is a slightly stronger result than that of \cite{Yuan}, although both are satisfactory for our purposes.
\medskip
\noindent In \cite{Corlette}, Corlette obtains a finiteness result for compact immersed submanifolds
of a given Riemannian manifold, recalling the finiteness result of Cheeger. As indicated by Corlette, this
finiteness result suggests an underlying compactness result. We introduce the following definition:
\proclaim{Definition \nextprocno}
\noindent Let $(M,g)$ be a Riemannian manifold. Let $X=(Y,i)$ be an immersed submanifold in $M$. Let $\nabla^i$ be the Levi-Civita covariant derivative generated over $Y$ by
the immersion $i$ into $(M,g)$. Let $A(X)$ be the Weingarten operator of $X$. For all $k\geqslant 2$, we define $A_k(X)$ using the following recurrence relation:
$$\matrix
A_2(X) \hfill&= A(X), \hfill\cr
A_k(X) \hfill&= \nabla^i A_{k-1}(X)\ \forall k\geqslant 3.\hfill\cr
\endmatrix$$
\endproclaim
\noindent If we denote by $\|\cdot\|_\Omega$ the $C^0$ norm over the set $\Omega$, then the result of 
Corlette may be stated as follows:
\proclaim{Theorem \nextprocno\ {\bf [Arzela, Ascoli, Corlette]}}
\noindent Let $(M_n,g_n,p_n)_{n\in\Bbb{N}}$ be a sequence of complete pointed manifolds which converges towards the pointed manifold $(M_0,g_0,p_0)$ in the
Cheeger/Gromov topology. For all $n\in\Bbb{N}$, let $\Sigma_n=(S_n,i_n)$ be a complete immersed submanifold in $M_n$, and suppose that there exists $q_n\in S_n$ such that
$i_n(q_n)=p_n$. Let $K\in\Bbb{N}$ be a positive integer. Suppose that, for all $R\in\Bbb{R}^+$ there exists
$B\in (0,\infty)$ and $N\in\Bbb{N}$ such that:
$$
n\geqslant N\Rightarrow \|A_k(\Sigma)\|_{B_R(q_n)}\leqslant B\qquad \forall k\leqslant K+1 .
$$
\noindent Then, there exists a complete pointed immersed submanifold $(\Sigma_0,q_0)=(S_0,i_0,q_0)$ in $M_0$ of type
$C^{K,\alpha}$ for all $\alpha\in(0,1)$ such that:
\medskip
\myitem{(i)} $i_0(q_0)=p_0$, and
\medskip
\myitem{(ii)} after extraction of a subsequence, $(\Sigma_n,p_n)_{n\in\Bbb{N}}$ converges to $(\Sigma_0,p_0)$ in the
$C^{K,\alpha}$\break Cheeger/Gromov topology for all $\alpha\in(0,1)$.
\endproclaim
\proclabel{LemmeArzelaAscoliCorlette}
\proof A proof may be found in \cite{SmiF}.\qed
\medskip
\noindent We also require the $\lambda$-maximum lemma. Let $X$ be a topological space. Let $f:X\rightarrow\Bbb{R}$ be a real valued function. We say that $f$ is {\emph locally bounded\/} if and only if for all $P\in X$ there exists $B\in\Bbb{R}^+$ and a neighbourhood $\Omega$ of $P$ in $X$ such that:
$$
Q\in\Omega\qquad\Rightarrow\qquad \left|f(Q)\right| \leqslant B.
$$
\goodbreak
\proclaim{Lemma \nextprocno\ {\sl $\lambda$-maximum lemma}}
\noindent Let $X$ be a metric space and let $P$ be a point in $X$. Suppose that there exists $\delta>0$ such that the closed $\delta$-ball about $P$ in $X$ is complete.
There exists $\lambda\in\Bbb{R}^+$, which only depends on $\delta$ such that, for every locally bounded function $f:X\rightarrow\Bbb{R}^+$ such that $f(P)\geqslant 1$, there
exists $Q\in B(P,\delta)$ such that $f(Q)\geqslant f(P)$ and:
$$
Q'\in B(Q,\epsilon)\Rightarrow f(Q')\leqslant \lambda f(Q),
$$
\noindent where $\epsilon$ is given by the relation:
$$
\epsilon = \lambda^{-1}f(Q)^{-1/2}.
$$
\endproclaim
\proclabel{LemmeDuLambdaMaximum}
\proof It suffices to assume the contrary and thus obtain a point in $X$ at which $f$ is arbitrarily large, which is absurd. A detailed proof may be found, for example, in \cite{SmiE}.\qed
\newsubhead{The Compactness Theorem}
\noindent Bearing in mind Theorem \procref{LemmeArzelaAscoliCorlette}, Theorem \procref{ThmCompactness} is
equivalent to the following result:
\proclaim{Lemma \nextprocno}
\noindent Let $M$ be a complete Riemannian manifold of dimension $2n+1$. Let $\theta\in (-n\pi/2,n\pi/2)$ be an angle.
Let $(W,m,\omega,\Omega_\theta)$ be a positive special Legendrian structure over $M$.
Let $\pi:\Cal{SL}^+_\theta(M)\rightarrow M$ be the canonical projection.
\medskip
\noindent For every compact subset $K$ of $M$, for every $R\in\Bbb{R}^+$ and for every $N\in\Bbb{N}$, there exists $B\in\Bbb{R}^+$ such that if $(\hat{\Sigma},p)\in\pi^{-1}(K)$ then:
$$
\|A_k(\hat{\Sigma})\|_{B_R(p)} \leqslant B\qquad \forall k\leqslant N.
$$
\endproclaim
\proclabel{LemmaCompactness}
\proof We assume the contrary and obtain a contradiction.
\medskip
{\noindent\bf First step : We begin by constructing a sequence of pointed submanifolds.}
\medskip
\noindent We may assume that there exists a compact subset $K$ of $M$, a sequence $(\Sigma_m,q_m)_\minn=(S_m,i_m,q_m)_\minn$
of pointed immersed submanifolds in $\Cal{SL}^+_\theta(M)$, and $k\geqslant 2$ such that $q_m\in K$ for all $m$, and:
$$
(\|A_k(\Sigma_m)(q_m)\|)_\minn\rightarrow\infty.
$$
\noindent For all $m$, we define $\Cal{A}_{k_m}:S_m\rightarrow\Bbb{R}$ by:
$$
\Cal{A}_{k,m} = \sum_{i=2}^k \|A_i(\hat{\Sigma}_m)\|^{\frac{1}{i-1}}.
$$
\noindent For all $m$, for all $\rho$ and for all $p\in S_m$, let $B(p,\rho)$ be the ball of radius $\rho$ about $p$ in $S_m$ with respect the the metric $i^*g$. Let $\lambda>1$ be a real number greater than $1$.
For all $m$, by replacing $q_m$ by a point $p_m$ obtained using the $\lambda$-maximum lemma (Lemma 
\procref{LemmeDuLambdaMaximum}), and by denoting $B_m=\Cal{A}_{k,m}(p_m)$, we may suppose that:
$$
q\in B\left(p_m,\frac{1}{\lambda B_m^{1/2}}\right) \Rightarrow \Cal{A}_{k,m}(q) \leqslant\lambda B_m.
$$
{\noindent\bf Second step : We simplify the problem by transforming the contact structure of $M$ onto the canonical contact structure of $\Bbb{R}\times\Bbb{R}^n\times\Bbb{R}^n$.}
\medskip
\noindent Since $K$ is compact, we may suppose that there exists $p_0$ to which $(p_m)_\minn$ converges. Let $\omega$ be
the canonical symplectic form over $\Bbb{R}^n\times\Bbb{R}^n$. Let $\beta$ be a primitive of $\omega$. We define the
contact form $\alpha$ over $\Bbb{R}\times\Bbb{R}^n\times\Bbb{R}^n$ by:
$$
\alpha = dt - \beta.
$$
\noindent By Darboux's theorem for families (see, for example, appendix $C$ of \cite{SmiE}) we may assume that there exists $\epsilon>0$ and, for all $m$, a neighbourhood
$U_m$ of $p_m$ in $M$ and a diffeomorphism $\varphi_m:(B_\epsilon(0),0)\rightarrow (U_m,p_m)$ such that
$(\varphi_m)_*\alpha$ is colinear with the contact structure over $U_m$. Moreover, we may assume that
$(\varphi_m)_\minn$ converges to $\varphi_0$ in the $C^\infty_{\oploc}$ topology.
\medskip
\noindent We denote by $\Bbb{R}^n_0$ the zero section over $\Bbb{R}^n$ in the bundle
$\Bbb{R}\times\Bbb{R}^n\times\Bbb{R}^n$. Since we may freely choose $\beta$, we may suppose that, for all $m$, the
application $T_0\varphi_m$ sends $\Bbb{R}^n_0$ onto $T_{p_m}\Sigma_m$.
\medskip
{\noindent\bf Third step : We rescale by dilating the metric over $M$ by a constant factor which tends to infinity as $m$ tends to infinity. The sequence of pointed manifolds thus obtained converges to a real vector space.}
\medskip
\noindent For all $m$, we define the metric $g_m$ over $M$ by:
$$
g_m = B_m^2g.
$$
\noindent For all $m$, we define $\Delta_m$, an endomorphism of $\Bbb{R}\times\Bbb{R}^n\times\Bbb{R}^n$, by:
$$
\Delta_m(t,q,p) = (B_m t, B_m q, B_m p).
$$
\noindent For all $m$, we define $\phi_m:(B_{B_m\epsilon}(0),0)\rightarrow(M,p_m)$ by:
$$
\phi_m = \varphi_m\circ\Delta_m^{-1}.
$$
\noindent Let $R$ be a positive real number. Since $(B_m)_\minn$ tends to infinity, there exists $M\in\Bbb{N}$ such that,
for $m\geqslant M$:
$$
B_R(0) \subseteq B_{B_m\epsilon}(0).
$$
\noindent Consequently, for $m\geqslant M$, $\phi_m$ is defined over $B_R(0)$ and is a diffeomorphism onto its image. Next:
$$\matrix
\phi_m^*g_m(t,q,p) \hfill&= ((\Delta_m^{-1})^*\varphi_m^* g_m)(t,q,p) \hfill\cr
&= (\varphi_m^*g)(t/B_m, q/B_m, p/B_m). \hfill\cr
\endmatrix$$
\noindent Since $(\varphi_m^*g)_\minn$ converges in the $C^\infty_\oploc$ topology towards $\varphi_0^* g$, the
sequence of metrics $(\phi_m^*g_m)_\minn$ converges to a constant metric $g_0$ over
$\Bbb{R}\times\Bbb{R}^n\times\Bbb{R}^n$. For the same reason, there exists a constant distribution $W_0$ in
$T\Bbb{R}\times\Bbb{R}^n\times\Bbb{R}^n$, a constant complex structure $J_0$ over $W_0$, a constant symplectic form
$\omega_0$ over $W_0$, a constant Minkowski metric $m_0$ over $W_0$, and a constant special Lagrangian form
$\Omega_0$ over $W_0$ such that $(\phi_m^*W)_\minn$, $(\phi_m^*J)_\minn$, $(B_m^2\phi_m^*\omega)_\minn$,
$(B_m^2\phi_m^*m)_\minn$ and $(B_m^n\phi_m^*\Omega)_\minn$ converge towards $W_0$, $J_0$, $\omega_0$, $m_0$ and
$\Omega_0$ respectively in the $C^\infty_\oploc$ topology.
\medskip
\goodbreak{\noindent\bf Fourth step : We transform the pointed submanifolds living in $M$ into submanifolds of the dilated
ambiant manifolds, and we show that the sequence of pointed submanifolds thus obtained converges towards an affine
subspace of a real vector space in the Cheeger/Gromov topology.}
\medskip
\noindent For all $m$, and for all $\rho\in (0,\infty)$ we denote by $B(p_m,\rho)$ the ball of radius $\rho$ about
$p_m$ in the manifold $S_m$ with respect to the metric $i^*g$ generated by the immersion $i$ into the
Riemannian manifold $(M, g)$. For sufficiently large $m$, we have:.
$$
B\left(p_m,\frac{1}{\lambda B_m^{1/2}}\right) \subseteq U_m.
$$
\noindent Consequently, for sufficiently large $m$, we define the pointed immersed submanifold
$\tilde{\Sigma}_m=(\tilde{S}_m,\tilde{\mathi}_m)$ which is contained in $\Bbb{R}\times\Bbb{R}^n\times\Bbb{R}^n$ by:
$$
(\tilde{S}_m,\tilde{\mathi}_m) = \left(B\left(p_m,\frac{1}{\lambda B_m^{1/2}}\right), \Delta_m\circ\varphi_m^{-1}\circ i_m\right).
$$
\noindent For all $m$, we view $\tilde{\Sigma}_m$ as an immersed submanifold in $(B_{B_m\epsilon}(0),0,\varphi_m^*g_m)$.
Thus, if we define $\tilde{\Cal{A}}_{k,m}$ for $\tilde{\Sigma}_m$ in the same way as $\Cal{A}_{k,m}$ for $\Sigma_m$, we obtain:
$$
\tilde{\Cal{A}}_{k,m}(p_m) = 1,
$$
\noindent and, for all $q\in S_m$:
$$
\tilde{\Cal{A}}_{k,m}(q) \leqslant \lambda.
$$
\noindent By the Arzela-Ascoli-Corlette Theorem (Theorem \procref{LemmeArzelaAscoliCorlette}), we may assume that
there exists a complete pointed immersed submanifold $\tilde{\Sigma}_0 = (\tilde{S}_0,\tilde{\mathi}_0)$ in
$\Bbb{R}\times\Bbb{R}^n\times\Bbb{R}^n$ of type $C^{k-1,\alpha}$ for all $\alpha$ which passes by the origin such that
$(\tilde{\Sigma}_m,p_m)_\minn$ tends towards $(\tilde{\Sigma}_0,p_0)$ in the $C^{k-1,\alpha}$-Cheeger/Gromov
topology for all $\alpha$.
\medskip
\noindent By using elliptic regularity along Theorem \procref{ThmYuanPlanes}, we obtain the
following result:
\proclaim{Lemma \nextprocno}
\noindent $\tilde{\Sigma}_0$ is a linear subspace of $\left\{0\right\}\times\Bbb{R}^n\times\Bbb{R}^n$ and
$\suite{\tilde{\Sigma}}{m}$ converges towards $\tilde{\Sigma}_0$ in the $C^\infty$ Cheeger/Gromov topology.
\endproclaim
\proclabel{LemmeRegulariteElliptique}
\proof $T\tilde{\Sigma}_0$ is contained in $W_0$. Since $W_0$ is a constant distribution over
$\Bbb{R}\times\Bbb{R}^n\times\Bbb{R}^n$, and since $\tilde{\Sigma}_0$ passes by the origin, it follows that
$\tilde{\Sigma}_0$ is contained in a linear subspace $X$ of $\Bbb{R}\times\Bbb{R}^n\times\Bbb{R}^n$ of dimension
$2n$ which is parallel to $W_0$. In fact, by construction of $W_0$, we obtain:
$$
X = \left\{0\right\}\times\Bbb{R}^n\times\Bbb{R}^n.
$$
\noindent The restriction of $(\omega_0,m_0,\Omega_0)$ defines a positive special Lagrangian structure over $X$. Since
$\tilde{\Sigma}_m$ is a positive special Lagrangian submanifold for all $m$, it follows that $\tilde{\Sigma}_0$ is also a positive
special Lagrangian submanifold. Since $\tilde{\Sigma}_0$ is complete and of type $C^{1,\alpha}$, it follows by
Jost and Xin's theorem (Theorem \procref{ThmYuanPlanes}) that $\tilde{\Sigma}_0$ is an affine subspace of $X$.
\medskip
\noindent We denote by $\Bbb{R}^n_0$ the zero section in the vector bundle $\Bbb{R}\times\Bbb{R}^n\times\Bbb{R}^n$ over
$\Bbb{R}^n$. By construction $T_0\tilde{\Sigma}_0$ coincides with $\Bbb{R}^n_0$ and it follows that $\tilde{\Sigma}_0$
itself coincides with $\Bbb{R}^n_0$. Since the sequence $(\tilde{\Sigma}_m)_\minn$ converges towards $\tilde{\Sigma}$
in the $C^{k-1,\alpha}$ Cheeger/Gromov topology, we may assume that there exists $\suite{\rho}{m}\in ]0,\infty[$ such that
$(\rho_m)_\minn\rightarrow\infty$ and, for all $m$, a function
$(t_m,q_m):B_{\rho_m}(0)\subseteq\Bbb{R}^n\rightarrow\Bbb{R}\times\Bbb{R}^n$ and a neighbourhood $\Omega_m$ about $p_m$ in
$S_m$ such that the image of $\Omega_m$ by $\tilde{\mathi}_m$ coincides with the graph of $(t_m,q_m)$ over $B_{\rho_m}(0)$.
The sequence of functions $(t_m,q_m)_\minn$ converges to zero in the $C^{k-1,\alpha}_\oploc$ topology. It suffices for our
needs that this sequence converges in the $C^{1,\alpha}_\oploc$ topology.
\medskip
\noindent For all $n\in\Bbb{N}$, we define the function $\Cal{O}_m(p,t,q,A)$ for $p,q\in\Bbb{R}^n$, $t\in\Bbb{R}$ and
$A\in\Bbb{R}^{n^2}$ with $\|p\|,\|q\|,\left|t\right|<\epsilon B_m/3$ by:
$$\matrix
\Cal{O}_m(p,t,q,A) \hfill&= (\Omega_m)_{\left(\frac{t}{B_m},\frac{p}{B_m},\frac{q}{B_m}\right)}\left[
\left(\partial_1, -\beta_{(p,q)}(\partial_1,A\cdot\partial_1), A\cdot\partial_1\right), ..., \right.\hfill\cr
&\qquad \left.\left(\partial_m, -\beta_{(p,q)}(\partial_m,A\cdot\partial_m), A\cdot\partial_m\right)\right].\hfill\cr
\endmatrix$$
\noindent Since $\tilde{\Sigma}_m$ is Legendrian with respect to the form $\alpha = dt - \beta$, we obtain, for all $m$:
$$
(dt_m)_p(\partial_i) = \beta_{(p,q_m(p))}(\partial_i, Dq_m\cdot\partial_i).
$$
\noindent Since $\tilde{\Sigma}_m$ is special Legendrian with respect to the form $\Omega_m$, we obtain, for all $m$:
$$
\Cal{O}_m(p, t_m, q_m, Dq_m) = 0.
$$
\noindent For all $\theta\in\Bbb{R}$, we define $\Cal{O}_{0,\theta}$ over
$\Bbb{R}^{(n+1)^2}=\Bbb{R}\times\Bbb{R}^n\times\Bbb{R}^n\times\Bbb{R}^{n^2}$ by:
$$
\Cal{O}_{0,\theta}(p,t,q,A) = \opIm(e^{i\theta}\opDet(I + iA)).
$$
\noindent There exists $\theta\in\Bbb{R}$ such that the sequence of functions $\suite{\Cal{O}}{m}$ converges towards
$\Cal{O}_{0,\theta}$ in the $C^\infty_\oploc$ topology over $\Bbb{R}^{(m+1)^2}$. Without loss of generality,
we may assume that $\theta=0$ and we define $\Cal{O}_0=\Cal{O}_{0,0}$. For all $m$, by taking the derivative
of $\Cal{O}_m$, we obtain for all $i$ a relation of the following form for $q_{m,i}$:
$$
a_{m,i}^{jk}\partial_j\partial_k q_{m,i} + b_{m,i}^j\partial_j q_{m,i} + c_{m,i} q_{m,i} = d_{m,i},
$$
\noindent where $a_{m,i}^{jk}$, $b_{m,i}^j$, $c_{m,i}$ and $d_{m,i}$ are smooth functions of $(p,t_m,q_m,Dq_m)$. For
all $i$, we define $a_{0,i}$ by:
$$
a_{0,i}^{pq}(p,t,q,A) = ((I+A^2)^{-1})^{pq}.
$$
\noindent Since $(\Cal{O}_m)_\minn$ tends towards $\Cal{O}_0$ in the $C^\infty_{\oploc}$ topology, the sequences
$(b^j_{m,i})_\minn$, $(c_{m,i})_\minn$ and $(d_{m,i})_\minn$ tend towards zero in the $C^\infty_\oploc$ topology, and
$(a^{jk}_{m,i})_\minn$ tends towards $a^{jk}_{0,i}$ in the $C^\infty_{\oploc}$ topology. Since $(t_m,q_m)_\minn$ converges
towards zero in the $C^{1,\alpha}_\oploc$ topology, it follows that the functions $a$,$b$,$c$ et $d$ are locally uniformly
bounded in the $C^{0,\alpha}$ norm. Moreover, $a_0$ is locally bounded from below, and these relations are
consequently uniformly elliptic. The Schauder estimates (see, for example, \cite{GilbTrud}) permit us to
conclude that, for all
$i$, $(D^2q_{m,i})_{n\in\Bbb{N}}$ is locally uniformly bounded in the $C^{0,\alpha}$ norm. It follows that $\suite{q}{m}$
is locally uniformly bounded in the $C^{2,\alpha}$ norm. Since, for all $m$, $t_m$ may be obtained from $q_m$ by
integrating the form $\beta(p,q_m)$, it follows that $\suite{t}{m}$ is also locally uniformly bounded in the $C^{2,\alpha}$
norm. Working by induction, we may thus show that $(t_m,q_m)_\minn$ is locally uniformly bounded in the $C^{k,\alpha}$ norm
for all $k\in\Bbb{N}$.
\medskip
\noindent The classical Arzela-Ascoli theorem now permits us to show that every subsequence of $(t_m,q_m)_\minn$ contains a subsubsequence which converges in the $C^k_\oploc$ topology for all $k$. The limit is necessarily zero. It thus follows that $(t_m,q_m)_\minn$ converges towards zero in the $C^\infty_\oploc$ topology as $m$ tends to infinity, and the result now follows.\qed
\medskip
{\noindent\bf Final step : We use continuity in order to obtain a contradiction.}
\medskip
\noindent We have shown that $(\tilde{\Sigma}_m,p_m)_\minn$ converges towards $(\tilde{\Sigma}_0,0)$ in the
Cheeger/Gromov topology. Moreover $(\tilde{\Sigma}_0,0)$ is an affine plane. However, if we define
$\tilde{\Cal{A}}_{k,0}$ over $\tilde{\Sigma}_0$ as before, since $(\tilde{\Sigma}_m,p_m)$ converges towards
$(\tilde{\Sigma}_0,0)$ in particular in the $C^k$-Cheeger-Gromov topology, we obtain:
$$
\tilde{\Cal{A}}_{k,0}(0) = 1.
$$
\noindent This is absurd. We thus obtain the desired contradiction, and the result now follows.\qed
\medskip
\noindent We now obtain Theorem \procref{ThmCompactness} as an immediate corollary:
\medskip
{\bf\noindent Proof of Theorem \procref{ThmCompactness}:} It follows by the preceeding lemma and the 
Arzela-Ascoli-Corlette Theorem (Theorem \procref{LemmeArzelaAscoliCorlette}) that there exists a subsequence of
$(\Sigma_n,p_n)_\ninn$ which converges in the $C^k$-Cheeger/Gromov topology for all $k$. In the case of immersed submanifolds, this is equivalent to convergence in the $C^\infty$-Cheeger-Gromov topology, and the
result now follows.\qed
\newhead{Degenerate Limits}
\newsubhead{The Derivatives of the Weingarten Operator}
\def\matA #1#2{{\Bbb{A}^{#1}}_{#2}}%
\def\hatnabla{\hat{\nabla}}%
\def\christ #1#2#3{{\Gamma^{#1}}_{#2#3}}%
\noindent Let $M$ be a Riemannian manifold. Let $g$ denote the metric over $M$. Let $\rho\in (0,\infty)$ be a positive
real number. Let $S_\rho M$ be the $\rho$-sphere bundle over $M$. Let $\Sigma=(S,i)$ be a convex, oriented, immersed
hypersurface in $M$. Let $\hat{\Sigma}_\rho=(S,\hat{\mathi}_\rho)$ be the $\rho$-Gauss lifting of $\Sigma$ in
$S_\rho M$. Let $A$ be the Weingarten operator of $\Sigma$. We define $\Bbb{A}$ by $\Bbb{A}=\rho A$.
\medskip
\noindent Let $p$ be a point in $S$. We work locally in a chart about $p$ in $S$ chosen such that the basis $\partial_1, ..., \partial_n$ is orthonormal at the origin, and the matrix $\Bbb{A}$ is diagonal at the origin.
\medskip
\noindent The matrix $\matA ij$ is a $(1,1)$ tensor. In the sequel, for $T^\alpha_\beta$ an arbitrary tensor over $S$,
we denote by $T^\alpha_{\beta;i}$ the covariant derivative of $T^\alpha_\beta$ with respect to the Levi-Civita
connection of $i^*g$ in the direction $\partial_i$. Raising and lowering of indices will be carried out with respect
to the metric $g_{ij} = i^*g$ (i.e. the metric inherited from $M$ through the immersion $i$). Let $\lambda_1, ..., \lambda_n$ be the eigenvalues of the matrix $\Bbb{A}$ at the origin. Using reasoning inspired by Calabi \cite{CalabiB}, we establish various relations concerning the symmetry of the derivatives of $\Bbb{A}$.
\medskip
\noindent Breaking slightly with tradition, for an arbitrary curvature tensor $R$,
we denote:
$$
R_{ijkl} = \langle R_{\partial_i\partial_j}\partial_k,\partial_l\rangle.
$$
\proclaim{Lemma \nextprocno}
\noindent Let $R^M$ and $R^\Sigma$ be the Riemann curvature tensors of $M$ and $\Sigma$ respectively. Let 
$\eta$ be the exterior normal vector to $\Sigma$. Then:
$$\matrix
\Bbb{A}_{ik;j} - \Bbb{A}_{jk;i} \hfill&= -\rho R^M_{ij\eta k},\hfill\cr
\Bbb{A}_{ij;kl} - \Bbb{A}_{ij;lk} \hfill&= R^{\Sigma}_{ipkl}\matA pj - R^\Sigma_{pjkl}\matA pi.\hfill\cr
\endmatrix$$
\endproclaim
\proclabel{DALemmaSymmOfDA}
\proclabel{DALemmaSymmOfDDA}
\proof This follows trivially from the definition of curvature (see, for example, \cite{SmiE}).\qed
\proclaim{Lemma \nextprocno}
\noindent If $\Sigma$ is of constant $\rho$-special Lagrangian curvature, then at the origin:
$$\matrix
\sum_{p=1}^n\frac{1}{1+\lambda_p^2}\Bbb{A}_{pp;k} \hfill&= 0,\hfill\cr
\sum_{p=1}^n\frac{1}{1+\lambda_p^2}\Bbb{A}_{pp;kl} \hfill&= 
\sum_{p,q=1}^n\frac{\lambda_p+\lambda_q}{(1+\lambda_p^2)(1+\lambda_q^2)}\Bbb{A}_{pq;k}\Bbb{A}_{pq;l}.\hfill\cr
\endmatrix$$
\endproclaim
\proclabel{ConstSLLemmaDA}
\proof This follows by differentiating the relation $\opIm(e^{-i\theta}\opDet(I+i\Bbb{A}))=0$.\qed

\newsubhead{Covariant Derivatives and Laplacians}
\noindent Let $\nabla$ denote the Levi-Civita covariant derivative of $(S,i^*g)$. Let $\hatnabla$ denote the Levi-Civita
covariante derivative of $(S,\hat{\mathi}_\rho^*g)$. Let $\christ kij$ denote the Christophel symbol of this covariant
derivative with respect to $\nabla$:
\subheadlabel{CovariantDerivatives}
$$
\christ kij\partial_k = (\hatnabla -\nabla)_{\partial_i}\partial_j.
$$
\noindent We observe that:
$$
\hat{\mathi}_\rho^*g(\cdot,\cdot) = i^*g((I+\Bbb{A}^2)\cdot,\cdot).
$$
\noindent We denote by $\Bbb{B}^{ij}$ the dual metric to $\hat{\mathi}_\rho^*g$. Thus:
$$
\Bbb{B}^{ij}(I+\Bbb{A}^2)_{jk} = {\delta^i}_k.
$$
\proclaim{Lemma \nextprocno}
\noindent $\Gamma$ satisfies the following equation:
$$
\christ kij=\frac{\rho}{2}\matA pj\Bbb{B}^{kq}R^M_{qip\eta} + \frac{\rho}{2}\matA pi\Bbb{B}^{kq}R^M_{qjp\eta}
+\frac{1}{2}\Bbb{B}^{kl}\matA pl(\Bbb{A}_{pi;j}+\Bbb{A}_{pj;i}).
$$
\endproclaim
\proclabel{CDLapLemmaChristSymbol}
\proof Analogously to the Koszul formula, we have:
$$\matrix
2\langle (\hat{\nabla} - \nabla)_{\partial_i}\partial_j, (I+\Bbb{A}^2)\partial_k \rangle \hfill
&=\langle\partial_j, [\nabla_{\partial_i}(I+\Bbb{A}^2)]\partial_k \rangle -
\langle\partial_i, [\nabla_{\partial_k}(I+\Bbb{A}^2)]\partial_j \rangle \hfill\cr
&\qquad + \langle\partial_k, [\nabla_{\partial_j}(I+\Bbb{A}^2)]\partial_i \rangle. \hfill\cr
\endmatrix$$
\noindent Thus:
$$
2(I+\Bbb{A}^2)_{kp}{\Gamma^p}_{ij}={\Bbb{A}^p}_j(\Bbb{A}_{pk;i} - \Bbb{A}_{pi;k}) + {\Bbb{A}^p}_i(\Bbb{A}_{pk;j} - \Bbb{A}_{pj;k}) + {\Bbb{A}^p}_k(\Bbb{A}_{pi;j}+ \Bbb{A}_{pj;i}).
$$
\noindent The result now follows by Lemma \procref{DALemmaSymmOfDA}.\qed
\medskip
\noindent Let $v$ be a point in $S_\rho M$. For $X,Y,Z\in T_vS_\rho M$, we define $\Cal{R}^a(X,Y,Z)$ and
$\Cal{R}^b(X,Y,Z)$ by:
$$\matrix
\Cal{R}^a(X,Y,Z) \hfill&= R^M(v, \pi_H(X), \pi_V(Y), \pi_H(Z)),\hfill\cr
\Cal{R}^b(X,Y,Z) \hfill&= R^M(v, \pi_V(X), \pi_H(Y), \pi_H(Z)).\hfill\cr
\endmatrix$$
\noindent Let $P$ be an $n$ dimensional subspace of $W_v\subseteq TS_\rho M$. Let $B^{ij}$ be the metric on $P^*$ dual to that
inherited by $P$ from $W_v$ through the canonical immersion. Let $e_1, ..., e_n$ be a basis for $P$. We define the
vector $\xi(P)=\xi(P)^k e_k\in P$ by:
$$
\xi(P)^k = - B^{ij}B^{kq}\Cal{R}^a_{iqj} - B^{ij}B^{kq}\Cal{R}^b_{iqj}.
$$
\noindent $\xi(P)$ does not depend on the basis of $P$ chosen. We thus define the vector field
$\tilde{\xi}$ over $TS\cong T\hat{\Sigma}$ by:
$$
\tilde{\xi}(p) = -\hat{\mathi}_\rho^*\xi(T_p\hat{\Sigma}_p).
$$
\noindent In the sequel, we will denote $\tilde{\xi}$ by $\xi$. In terms of the basis
$\partial_1, ...,\partial_n$ of $TS$:
$$
\xi^q = -\rho\Bbb{B}^{ij}\Bbb{B}^{kq}\matA pj R^M_{\eta ipk} - \rho\Bbb{B}^{ij}\Bbb{B}^{kq}\matA pi R^M_{\eta pjk}.
$$
\noindent Since $\xi$ is defined in terms of $\hat{\mathi}_\rho$, it remains meaningful as long as $\hat{\mathi}_\rho$ converges, even if $i$ doesn't.
\proclaim{Lemma \nextprocno}
\noindent For every function $f$ over $S$:
$$
\hat{\Delta}f + df(\xi) = \Bbb{B}^{ij}f_{;ij}.
$$
\endproclaim
\proclabel{CDLapLemmaLaplacian}
\proof First:
$$
\hat{\Delta} f = \Bbb{B}^{ij}f_{;ij} - \Bbb{B}^{ij}{\Gamma^k}_{ij}f_{;k}.
$$
\noindent By Lemma \procref{CDLapLemmaChristSymbol} and the symmetry of $\Bbb{B}$:
$$
\Bbb{B}^{ij}{\Gamma^k}_{ij} = -\rho\Bbb{B}^{jq}\Bbb{B}^{kp}{\Bbb{A}^i}_qR^M_{\eta ipj} +
\frac{1}{2}\Bbb{B}^{ij}\Bbb{B}^{kq}{\Bbb{A}^p}_q(\Bbb{A}_{pi;j} + \Bbb{A}_{pj;i}).
$$
\noindent By Lemma \procref{DALemmaSymmOfDA}:
$$
\Bbb{A}_{pi;j} + \Bbb{A}_{pj;i} = \Bbb{A}_{ji;p} + \Bbb{A}_{ij;p} + \rho R^M_{\eta ijp} + \rho R^M_{\eta jip}.
$$
\noindent Using Lemma \procref{ConstSLLemmaDA} and the symmetry of $\Bbb{B}$ yields:
$$
\Bbb{B}^{ij}(\Bbb{A}_{pi;j} + \Bbb{A}_{pj;i}) = 2\rho\Bbb{B}^{ij}R^M_{\eta ijp}.
$$
\noindent Thus:
$$
\Bbb{B}^{ij}{\Gamma^k}_{ij} = -\rho\Bbb{B}^{ij}\Bbb{B}^{kq}{\Bbb{A}^p}_qR^M_{\eta ipj}
-\rho\Bbb{B}^{jq}\Bbb{B}^{kp}{\Bbb{A}^i}_qR^M_{\eta ipj} = \xi^k.
$$
\noindent The result now follows.\qed
\newsubhead{A Useful Function}
\noindent Let $HS_\rho M$ be the horizontal bundle of $S_\rho M$ associated to the Levi-Civita connection over $M$. Let $\pi_H$ be the projection of $W_\rho$ onto $HS_\rho M$. For all $\tau\in(0,\infty)$, we define
$f^\tau$ by:
\subheadlabel{UsefulFunction}
$$
f^\tau = \opDet(\pi_H)^{\tau/2} = \opDet(I+\Bbb{A}^2)^{-\tau}.
$$
\noindent For all $\tau$, the function $f^\tau$ is defined even when $\hat{\Sigma}$ is vertical. We aim to calculate the Laplacian of $f^\tau$. We first require:
\proclaim{Lemma \nextprocno}
\noindent Let $\eta$ be the exterior normal vector to $\Sigma$. Let $\nabla^M R^M$ be the covariant derivative of $R^M$ with respect to the Levi-Civita connection over $M$. Then:
$$
\Bbb{A}_{ij;kl} - \Bbb{A}_{kj;il} = -\rho(\nabla^MR^M)_{ik\eta jl} - R^M_{\eta k\eta j}\Bbb{A}_{li}- R^M_{i\eta\eta j}\Bbb{A}_{lk} - \matA plR^M_{ikpj}.
$$
\endproclaim
\proclabel{UsefulLemmaSymmOfDDA}
\proof This follows by differentiating the first relation in Lemma \procref{DALemmaSymmOfDA}.\qed
\medskip
\noindent For all $i$, we define $\hat{\partial}_{i,\rho}$ by:
$$
\hat{\partial}_{i,\rho} = (\hat{\mathi}_\rho)_*\partial_i = \left\{\partial_i,\rho A\partial_i\right\}.
$$
\noindent Let $\pi:TS_\rho M\rightarrow W_\rho$ be the orthogonal projection onto $W_\rho$. We denote by
$\nabla^{W_\rho}=\pi\circ\nabla^{S_\rho M}$ the covariant derivative of the distribution $W_\rho$ inherited
from the Levi-Civita connection of $S_\rho M$. We define $\hat{A}^\rho_{ijk}$ over $S$ by:
$$\matrix
\hat{A}^\rho_{ijk} \hfill&= \omega(\nabla^{W_\rho}_{\hat{\partial}_{i,\rho}}\hat{\partial}_{j,\rho},\hat{\partial}_{k,\rho})
+\frac{1}{2}R^M(\pi_H(\hat{\partial}_{i,\rho}),\pi_H(\hat{\partial}_{j,\rho}),v,\pi_H(\hat{\partial}_{k,\rho}))\hfill\cr
&+\frac{1}{2}R^M(v,\pi_V(\hat{\partial}_{j,\rho}),\pi_H(\hat{\partial}_{i,\rho}),\pi_V(\hat{\partial}_{k,\rho}))
+\frac{1}{2}R^M(v,\pi_V(\hat{\partial}_{i,\rho}),\pi_H(\hat{\partial}_{j,\rho}),\pi_V(\hat{\partial}_{k,\rho})).\hfill\cr
\endmatrix$$
\noindent In the sequel, we will refer to $\hat{A}^\rho$ as the {\emph adjusted Wiengarten operator\/} of $\hat{\Sigma}_\rho$. We observe that $\|\hat{A}^\rho\|$ is controlled by the norm of $R^M$ and the norm
of the second fundamental form of $\hat{\Sigma}_\rho$.
\noskipproclaim{Lemma \nextprocno}
$$
\hat{A}^\rho_{ijk} = \Bbb{A}_{jk;i}.
$$
\endnoskipproclaim
\proclabel{DALemmaRelationAandTildeA}
\proof This is a trivial calculation (see, for example, \cite{SmiE}).\qed
\noskipproclaim{Lemma \nextprocno}
\noindent Let $\hat{A}^\rho$ be the adjusted Weingarten operator of $\hat{\Sigma}_\rho$. If 
$h=\opLog(\opDet(I+\Bbb{A}^2))$, then:
$$\matrix
\Bbb{B}^{ij}h_{;ij} \hfill&=
\sum_{pq,r}\frac{2(1+\lambda_p\lambda_q)}{(1+\lambda_p^2)(1+\lambda_q^2)(1+\lambda_r^2)}\Bbb{A}_{pq;r}\Bbb{A}_{pq;r}\hfill\cr
&\qquad-\sum_{p,q}\frac{(\lambda_p-\lambda_q)^2}{(1+\lambda_p^2)(1+\lambda_q^2)}R^\Sigma_{pqpq}\hfill\cr
&\qquad-C(T_p\hat{\Sigma},M,\tilde{A}^\rho),\hfill\cr
\endmatrix$$
\noindent where $C$ is a continuous function of:
\medskip
\myitem{(i)} the affine subspace $T_p\hat{\Sigma}\subseteq W_p$,
\medskip
\myitem{(ii)} the tensor $\hat{A}^\rho$ defined over $T_p\hat{\Sigma}$, and
\medskip
\myitem{(iii)} the curvature of $M$ and its first derivative.
\endproclaim
\proclabel{LemmaLapOfLog}
\proof Differentiating $h$ yields:
$$\matrix
h_{;i} \hfill&= 2\opTr(\Bbb{A}(I+\Bbb{A}^2)^{-1}\Bbb{A}_{;i}),\hfill\cr
h_{;ij} \hfill&= 2\opTr(\Bbb{A}_{;j}(I+\Bbb{A}^2)^{-1}\Bbb{A}_{;i}) \hfill\cr
&\qquad -2\opTr(\Bbb{A}(I+\Bbb{A}^2)^{-1}[\Bbb{A}\Bbb{A}_{;j} + \Bbb{A}_{;j}\Bbb{A}](I+\Bbb{A}^2)^{-1}\Bbb{A}_{;i}) \hfill\cr
&\qquad +2\opTr(\Bbb{A}(I+\Bbb{A}^2)^{-1}\Bbb{A}_{;ij}).\hfill\cr
\endmatrix$$
\noindent Thus, at the origin:
$$\matrix
h_{;ij} \hfill&= \sum_{p,q}\frac{2}{1+\lambda_p^2}\Bbb{A}_{pq;i}\Bbb{A}_{pq;j} -
\sum_{p,q}\frac{2\lambda_p^2}{(1+\lambda_p^2)(1+\lambda_q^2)}\Bbb{A}_{pq;i}\Bbb{A}_{pq;j} \hfill\cr
&\qquad -\sum_{p,q}\frac{2\lambda_p\lambda_q}{(1+\lambda_p^2)(1+\lambda_q^2)}\Bbb{A}_{pq;i}\Bbb{A}_{pq;j} + \sum_p\frac{2\lambda_p}{1+\lambda_p^2}\Bbb{A}_{pp;ij}.\hfill\cr
\endmatrix$$
\noindent Using the symmetry of $\Bbb{A}$ with respect to $p$ and $q$, we obtain:
$$
h_{;ij} = \sum_{p,q}\frac{2-2\lambda_p\lambda_q}{(1+\lambda_p^2)(1+\lambda_q^2)}\Bbb{A}_{pq;i}\Bbb{A}_{pq;j} + \sum_p\frac{2\lambda_p}{1+\lambda_p^2}\Bbb{A}_{pp;ij}.
$$
\noindent We aim to eliminate the second derivatives in this expression.
$$\matrix
\sum_{p,q}\frac{\lambda_p}{1+\lambda_p^2}\frac{1}{1+\lambda_q^2}\Bbb{A}_{pp;qq} &= \sum_{p,q}\frac{\lambda_p}{(1+\lambda_p^2)(1+\lambda_q^2)}\Bbb{A}_{qq;pp}\hfill\cr
&\qquad + \sum_{p,q}\frac{\lambda_p}{(1+\lambda_p^2)(1+\lambda_q^2)}[\Bbb{A}_{qp;qp}-\Bbb{A}_{qq;pp}]\hfill\cr
&\qquad +\sum_{p,q}\frac{\lambda_p}{(1+\lambda_p^2)(1+\lambda_q^2)}[\Bbb{A}_{qp;pq}-\Bbb{A}_{qp;qp}]\hfill\cr
&\qquad +\sum_{p,q}\frac{\lambda_p}{(1+\lambda_p^2)(1+\lambda_q^2)}[\Bbb{A}_{pp;qq}-\Bbb{A}_{qp;pq}].\hfill\cr
\endmatrix$$
\noindent By Lemmata \procref{DALemmaSymmOfDDA} and \procref{UsefulLemmaSymmOfDDA}:
$$\matrix
\sum_{p,q}\frac{\lambda_p}{1+\lambda_p^2}\frac{1}{1+\lambda_q^2}\Bbb{A}_{pp;qq} \hfill&= \sum_{p,q}\frac{\lambda_p}{(1+\lambda_p^2)(1+\lambda_q^2)}\Bbb{A}_{qq;pp}\hfill\cr
&\qquad+\sum_{p,q}\frac{\lambda_p(\lambda_q-\lambda_p)}{(1+\lambda_p^2)(1+\lambda_q^2)}R^\Sigma_{pqpq}+C_1(T_p\hat{\Sigma},M),\hfill\cr
\endmatrix$$
\noindent where $C_1$ depends continuously on $T_p\hat{\Sigma}$, the curvature of $M$ and its first derivative. The last relation in Lemma \procref{ConstSLLemmaDA}
now yields:
$$\matrix
\sum_{p,q}\frac{\lambda_p}{1+\lambda_p^2}\frac{1}{1+\lambda_q^2}\Bbb{A}_{pp;qq} \hfill&=
\sum_{p,q,r}\frac{(\lambda_p+\lambda_q)\lambda_r}{(1+\lambda_p^2)(1+\lambda_q^2)(1+\lambda_r^2)}\Bbb{A}_{pq;r}\Bbb{A}_{pq;r}\hfill\cr
&\qquad + \sum_{p,q}\frac{\lambda_p(\lambda_q-\lambda_p)}{(1+\lambda_p^2)(1+\lambda_q^2)}R^\Sigma_{pqpq}+C_1(T_p\hat{\Sigma},M).\hfill\cr
\endmatrix$$
\noindent Thus:
$$\matrix
\Bbb{B}^{ij} h_{;ij} \hfill&= \sum_p\frac{1}{1+\lambda_p^2}h_{;pp} \hfill\cr
&= \sum_{p,q,r}\frac{2-2\lambda_p\lambda_q}{(1+\lambda_p^2)(1+\lambda_q^2)(1+\lambda_r^2)}\Bbb{A}_{pq;r}\Bbb{A}_{pq;r} +
2\sum_{p,q}\frac{\lambda_p}{(1+\lambda_p^2)(1+\lambda_q^2)}\Bbb{A}_{pp;qq} \hfill\cr
&= \sum_{p,q,r}\frac{2-2\lambda_p\lambda_q + 2\lambda_p\lambda_r + 2\lambda_q\lambda_r}{(1+\lambda_p^2)(1+\lambda_q^2)(1+\lambda_r^2)}\Bbb{A}_{pq;r}\Bbb{A}_{pq;r} \hfill\cr
&\qquad +2\sum_{p,q}\frac{\lambda_p(\lambda_q-\lambda_p)}{(1+\lambda_p^2)(1+\lambda_q^2)}R^\Sigma_{pqpq}+C_1(T_p\hat{\Sigma},M).\hfill\cr
\endmatrix$$
\noindent By the symmetry of $\Bbb{A}_{ij}$ and Lemma \procref{DALemmaSymmOfDA}:
$$\matrix
\sum_{p,q,r}\frac{2\lambda_p\lambda_r + 2 \lambda_q\lambda_r}{(1+\lambda_p^2)(1+\lambda_q^2)(1+\lambda_r^2)}\Bbb{A}_{pq;r}\Bbb{A}_{pq;r}\hfill&=
\sum_{p,q,r}\frac{4\lambda_q\lambda_r}{(1+\lambda_p^2)(1+\lambda_q^2)(1+\lambda_r^2)}\Bbb{A}_{pq;r}^2\hfill\cr
&=\sum_{p,q,r}\frac{4\lambda_q\lambda_r}{(1+\lambda_p^2)(1+\lambda_q^2)(1+\lambda_r^2)}(\Bbb{A}_{rq;p} - \rho 
R^M_{pr\eta q})^2.\hfill\cr
\endmatrix$$
\noindent Since:
$$
\lambda_q\lambda_rR^M_{pr\eta q} = R^M(\pi_H(\hat{\partial}_{i,\rho}),\pi_V(\hat{\partial}_{r,\rho}),v,\pi_V(\hat{\partial}_{q,\rho})),
$$
\noindent we obtain:
$$\matrix
\Bbb{B}^{ij} h_{;ij} \hfill&= \sum_{p,q,r}\frac{2+2\lambda_p\lambda_q}{(1+\lambda_p^2)(1+\lambda_q^2)(1+\lambda_r^2)}\Bbb{A}_{pq;r}\Bbb{A}_{pq;r} \hfill\cr
&\qquad +2\sum_{p,q}\frac{\lambda_p(\lambda_q-\lambda_p)}{(1+\lambda_p^2)(1+\lambda_q^2)}R^\Sigma_{pqpq}+C(T_p\hat{\Sigma},M,\hat{A}^\rho).\hfill\cr
\endmatrix$$
\noindent The result now follows by the antisymmetry of $R^\Sigma$.\qed
\medskip
\noindent We now obtain the following result concerning $f^\tau$ and $\hat{A}^\rho$:
\proclaim{Lemma \nextprocno}
\noindent With the same notation as in Lemma \procref{LemmaLapOfLog}, if $\tau\leqslant 1/2n$ then:
$$
\frac{1}{\tau f^\tau}\Bbb{B}^{ij}f^\tau_{;ij} \leqslant
\sum_{p,q}\frac{(\lambda_p-\lambda_q)^2}{(1+\lambda_p^2)(1+\lambda_q^2)}R^\Sigma_{pqpq}
+C(T_p\hat{\Sigma},M,\hat{A}^\rho).
$$
\endproclaim
\proclabel{UsefulLemmaLaplaceF}
\proof Trivially:
$$
\frac{1}{\tau f^\tau}\Bbb{B}^{ij}f_{;ij} = \frac{1}{\tau(f^\tau)^2}\Bbb{B}^{ij}f^\tau_{;i}f^\tau_{;j} - \Bbb{B}^{ij}h_{;ij}.
$$
\noindent The derivative of $f^\tau$ satisfies:
$$
f^\tau_{;i} = -\tau f^\tau\sum_{j=1}^n \frac{\lambda_j}{(1+\lambda_j^2)}\Bbb{A}_{jj;i}.
$$
\noindent Thus:
$$
\frac{1}{\tau(f^\tau)^2}\Bbb{B}^{ij}f^\tau_{;i}f^\tau_{;j}=
\tau\sum_{p,q,r}\frac{1}{(1+\lambda_r^2)}\frac{4\lambda_p\lambda_q}{(1+\lambda_p^2)(1+\lambda_q^2)}
\Bbb{A}_{pp;r}\Bbb{A}_{qq;r}.
$$
\noindent However, for all $r$, using the Cauchy-Schwarz inequality, we obtain:
$$\matrix
\sum_{p,q}\frac{4\lambda_p\lambda_q}{(1+\lambda_p^2)(1+\lambda_q^2)}\Bbb{A}_{pp;r}\Bbb{A}_{qq;r}
\hfill&= (\sum_p \frac{2\lambda_p}{(1+\lambda_p^2)}\Bbb{A}_{pp;r})^2\hfill\cr
&\leqslant 4n\sum_p \frac{\lambda_p^2}{(1+\lambda_p^2)^2}\Bbb{A}_{pp;r}^2\hfill\cr
&\leqslant 4n\sum_{p,q}\frac{\lambda_p\lambda_q}{(1+\lambda_p^2)(1+\lambda_q^2)}\Bbb{A}_{pq;r}\Bbb{A}_{pq;r}.\hfill\cr
\endmatrix$$
\noindent Consequently:
$$
\frac{1}{\tau(f^\tau)^2}B^{ij}f^\tau_{;i}f^\tau_{;j}\leqslant
2n\tau\sum_{p,q,r}\frac{2(1+\lambda_p\lambda_q)}{(1+\lambda_p^2)(1+\lambda_q^2)(1+\lambda_r^2)}\Bbb{A}_{pq;r}\Bbb{A}_{pq;r}.
$$
\noindent Thus, if $2n\tau\leqslant 1$, the desired result follows from Lemmata \procref{DALemmaRelationAandTildeA} and \procref{LemmaLapOfLog}.\qed
\medskip
\noindent To summarise, we obtain:
\proclaim{Lemma \nextprocno}
\noindent If $\tau\leqslant 1/2n$, then:
$$
\hat{\Delta} f^\tau + df^\tau(\xi) - Cf^\tau \leqslant 0.
$$
\noindent where $C>0$ is a constant depending continuously on the local geometries of $\hat{\Sigma}$ and $M$.
\endproclaim
\proclabel{Subharmonic}
\proof By definition of $\Bbb{A}$:
$$
R^\Sigma_{pqpq} = R^M_{pqpq} - \rho^{-2}\lambda_p\lambda_q(1-\delta_{pq}).
$$
\noindent Since $\lambda_i\geqslant 0$ for all $i$:
$$
\sum_{p,q}\frac{(\lambda_p-\lambda_q)^2}{(1+\lambda_p^2)(1+\lambda_q^2)}R^\Sigma_{pqpq} \leqslant \Sigma_{p,q}\frac{(\lambda_p-\lambda_q)^2}{(1+\lambda_p^2)(1+\lambda_q^2)}R^M_{pqpq}.
$$
\noindent As before, the term on the right hand side depends continuously on the curvature $M$ and may thus be incorporated into $C$. Thus, by Lemmata
\procref{CDLapLemmaLaplacian} and \procref{UsefulLemmaLaplaceF}:
$$
\hat{\Delta}f^\tau + df^\tau(\xi) - C(T_p\hat{\Sigma},M,\tilde{A}_p)f^\tau\leqslant 0.
$$
\noindent Since $f^\tau$ is positive, the result now follows by taking a constant upper bound for $C$.\qed
\newsubhead{Degeneration of Submanifolds}
\noindent We are now in a position to describe submanifolds in $\partial\hat{\Cal{F}_{\rho,\theta}(M)}(M)$. We require the following definition:
\proclaim{Definition \nextprocno}
\noindent Let $M$ be a manifold of dimension $(n+1)$. Let $S_\rho M$ be the $\rho$-sphere bundle over $M$. Let
$\pi:S_\rho M\rightarrow M$ be the canonical projection. Let $\hat{\Sigma}=(S,\hat{\mathi})$ be an $n$ dimensional immersed
submanifold of $S_\rho M$. Let $p$ be an arbitrary point of $S$. We say that $\hat{\Sigma}$ is {\emph vertical} at $p$ if
and only if the kernel of $T_{i(p)}\pi$ has non-trivial intersection with $T_p\hat{\Sigma}=T_p\hat{\mathi}\cdot T_pS$.
For $k\in\Bbb{N}$, we say that $\hat{\Sigma}$ is {\emph vertical of order $k$\/} at $p$ if and only if:
$$
\opDim(\opKer(T_{i(p)}\pi)\minter T_p\hat{\Sigma}) = k.
$$
\noindent We say that $\Sigma$ is {\emph vertical\/} if and only if it is vertical at every point $p$ in $S$.
\endproclaim
\proclaim{Lemma \nextprocno}
\noindent If $(\Sigma,p)$ is a pointed immersed submanifold in $\partial\hat{\Cal{F}}^+_{\rho,\theta}(M)$, then $\Sigma$ is vertical.
\endproclaim
\proclabel{LemmaBoundaryVerticalForThetaNotInteger}
\proof Let $(\hat{\Sigma}_0,p_0)=(S_0,\hat{\mathi}_0,p_0)$ be a pointed, immersed submanifold in $\partial\hat{\Cal{F}}_{\rho,\theta}(M)$.
Let $(\Sigma_m,p_m)_{m\in\Bbb{N}}=(S_m,i_m,p_m)_{m\in\Bbb{N}}$ be a sequence of convex, pointed, immersed submanifolds in
$\Cal{F}_{\rho,\theta}(M)$ such that if, for all $n$, $(\hat{\Sigma}_m,p_m)$ is the $\rho$-Gauss lifting of
$(\Sigma_m,p_m)$, then $(\hat{\Sigma}_m,p_m)_{m\in\Bbb{N}}$ converges to $(\Sigma_0,p_0)$ in the Cheeger/Gromov topology.
\medskip
\noindent For all $m$, let $A_m$ be the Weingarten operator of $\Sigma_m$, and let $\lambda_{m,1},...,\lambda_{m,n}\geqslant 0$ be
the eigenvalues of $\rho A_m(p_m)$. We define the function $f_m:S_m\rightarrow\Bbb{R}$ by:
$$
f_m = \opDet(I+\rho^2 A_m^2)^{-1/2n}.
$$
\noindent We define $f_0$ over $S_0$ by $f_0=f^{1/2n}(\Sigma_0)$. By Lemma \procref{Subharmonic}, for all $m$:
$$
\hat{\Delta}f_m + df_m(\xi_m) - Cf_m \leqslant 0.
$$
\noindent Since $(\hat{\Sigma}_m)_{m\in\Bbb{N}}$ converges, $C$ may be chosen the same for all $m$. By the maximum principle (see Lemma \procref{MaximumPrinciple}), if 
$f_0$ vanishes at a single point, then it vanishes everywhere. However, $\Sigma_0$ is vertical if and only if $f_0=0$ and the result follows.\qed
\newsubhead{The Geometry of Curtain Submanifolds}
\noindent We now prove Theorem \procref{ThmCompactnessOfHypersurfaces}. We recall the following result concerning vertical submanifolds:
\proclaim{Lemma \nextprocno}
\noindent Let $M$ be a Riemannian manifold. Let $\rho$ be a positive real number. Let $\hat{\Sigma}=(S,\hat{\mathi})$ be an
immersed submanifold in $S_\rho M$. Let $\pi:S_\rho M\rightarrow M$ be the canonical projection. Let $p$ be an arbitrary point in $S$.
Suppose that there exists a neighbourhood $\Omega$ of $p$ in $S$ such that $\opDim(T\hat{\Sigma}\minter VS_\rho M)$ is constant over
$\Omega$. Define $k$ by $k=\opDim(T\hat{\Sigma}\minter VS_\rho M)$. By restricting $\Omega$ if necessary, we may assume that
there exists $\epsilon>0$, a submanifold $\Sigma\subseteq M$ embedded into $M$ and a diffeomorphism
$\Phi:\Sigma\times(-\epsilon,\epsilon)^k\rightarrow\Omega$ such that, if we denote by
$\pi_1:\Sigma\times(-\epsilon,\epsilon)^k\rightarrow\Sigma$ the projection onto the first factor, then:
$$
\pi\circ\Phi = \pi_1.
$$
\endproclaim
\proclabel{GCSLemmaDecomposition}
\noindent We now obtain the following result:
\proclaim{Lemma \nextprocno}
\noindent Let $M$ be a Riemannian manifold. Let $\rho$ be a positive real number. Let $\theta\in[(n-1)\pi/2,n\pi/2[$ be an angle. Let $(\hat{\Sigma},p)=(S,\hat{\mathi},p)$ 
be an immersed submanifold in $\partial\hat{\Cal{F}}_{\rho,\theta}(M)$. If $\hat{\Sigma}$ is vertical, then
$\hat{\Sigma}$ is a normal sphere bundle over a complete totally geodesic submanifold of $M$.
\endproclaim
\proclabel{LemmaLimitIsTotallyGeodesicSphereBundle}
\proof Let $(\Sigma_m,p_m)_{m\in\Bbb{N}}=(S_m,i_m,p_m)_{m\in\Bbb{N}}$ be a sequence of immersed
submanifolds in $\hat{\Cal{F}}^+_{\rho,\theta}(M)$ whose $\rho$-Gauss liftings converge to $(\hat{\Sigma},p)$ in the Cheeger/Gromov
topology.
\medskip
\noindent Let $q$ be a point in $S$. Suppose that $\opDim(T\hat{\Sigma}\minter VS_\rho M)$ is constant in a neighbourhood of
$q$. By Lemma \procref{GCSLemmaDecomposition}, there exists a submanifold $\Sigma$ of $M$, $k\in\Bbb{N}$, $\epsilon\in\Bbb{R}^+$, a
neighbourhood $\Omega$ of $q$ in $S$ and a diffeomorphism $\Phi:\Sigma\times(-\epsilon,\epsilon)^k\rightarrow\Omega$ such that:
$$
(\pi\circ\hat{\mathi})\circ\Phi = \pi_1.
$$
\noindent Let $(x,y)$ be a point in $\Sigma\times (-\epsilon,\epsilon)^k$. Let $X_1,...,X_n$ be a basis of
$T_x\Sigma$. For all $i$, we have:
$$
T\pi\cdot T(\hat{\mathi}\circ\Phi)\cdot X_i = X_i.
$$
\noindent Consequently, there exists $Y_1, ..., Y_m\in T_xM$ such that, for all $i$:
$$
T(\hat{\mathi}\circ\Phi)\cdot X_i = \left\{X_i,Y_i\right\}.
$$
\noindent Let $(q_m)_{m\in\Bbb{N}}$ be a sequence of points in $(S_m)_{m\in\Bbb{N}}$ which tends towards
$\Phi(x,y)$. For all $m$, let $(E_{1,m}, ...,E_{n,m})$ be vectors in $T_{q_m}S_m$ such that, for all $i$:
$$\matrix
&(T\hat{\mathi}_m\cdot E_{i,m})_{m\in\Bbb{N}}\hfill&\rightarrow \left\{X_i,Y_i\right\}\hfill\cr
\Rightarrow\hfill& (T(\pi\circ\hat{\mathi}_m)\cdot E_{i,m})_{m\in\Bbb{N}}\hfill&\rightarrow X_i.\hfill\cr
\endmatrix$$
\noindent For all $i$ and for all $m$, let $X_{i,m},Y_{i,m}\in T_{(\pi\circ\hat{\mathi}_m)(q_m)}M$ be such that:
$$
T\hat{\mathi}_m\cdot E_{i,m} = \left\{ X_{i,m},Y_{i,m}\right\}.
$$
\noindent For all $m$, and for all $i$, since $X_{i,m}$ is a tangent vector of the hypersurface of which $\hat{\Sigma}_m$
is the $\rho$-Gauss lifting, we have:
$$
\langle \hat{\mathi}_m(q_m),X_{i,m}\rangle = 0.
$$
\noindent Thus, by taking limits, we obtain:
$$
\langle \hat{\mathi}(\Phi(x,y)),X_i\rangle = 0\qquad\forall i.
$$
\noindent Consequently, if we denote by $N\Sigma$ the normal $\rho$-sphere bundle over $\Sigma$, we find that, for all
$q'\in\Omega$:
$$
\hat{\mathi}(q')\in N\Sigma.
$$
\noindent Since $N\Sigma$ is of dimension $n$, it follows that $\hat{\mathi}$ sends $\Omega$ diffeomorphically onto an open subset of $N\Sigma$. Since $N\Sigma$ is $\theta$-special
Lagrangian, it follows by the unique continuation principle (Lemma \procref{ThmUniqueContinuation}) that $\Omega$ may be extended to an open set of $S$ such that $\hat{\mathi}:\Omega\rightarrow N\Sigma$ is a covering map.
\medskip
\noindent For all $i$ and for all $m$, since $\hat{\Sigma}_m$ is the $\rho$-Gauss lifting of a convex hypersurface, we
have:
$$
\langle X_{i,m},Y_{i,m}\rangle \geqslant 0.
$$
\noindent Consequently, taking limits, we obtain, for all $i$:
$$
\langle X_i, Y_i\rangle \geqslant 0.
$$
\noindent Thus, since $\hat{\mathi}$ is a diff\'eomorphism, for all $\left\{X,Y\right\}\in TN\Sigma$:
$$
\langle X, Y\rangle \geqslant 0.
$$
\noindent Therefore, if $N$ is a normal vector field of $\Sigma$, then, for all $X\in T_\Sigma$:
$$
\langle X,\nabla_X N\rangle \geqslant 0.
$$
\noindent However, $(-N)$ is also a normal vector field over $\Sigma$, and thus, for all $X\in T\Sigma$:
$$\matrix
&\langle X,\nabla_X N\rangle \hfill&= -\langle X, \nabla_X(-N)\rangle \hfill&\leqslant 0\hfill\cr
\Rightarrow\hfill& \langle X,\nabla_X N\rangle \hfill&=0.\hfill&\cr
\endmatrix$$
\noindent Thus, by polarisation, for all $X,Y\in T\Sigma$ and for all $N$ normal to $\Sigma$:
$$
\langle X, \nabla_Y N\rangle = 0.
$$
\noindent It thus follows that $\Sigma$ is a totally geodesic submanifold. The result now follows by unique continuation.\qed
\medskip
{\bf\noindent Proof of Theorem \procref{ThmCompactnessOfHypersurfaces}:\ }Let $(\hat{\Sigma},p)$ be an immersed submanifold in $\partial\hat{\Cal{F}}_{\rho,\theta}(M)$.
By Lemma \procref{LemmaBoundaryVerticalForThetaNotInteger}, $\hat{\Sigma}$ is vertical. By Lemma
\procref{LemmaLimitIsTotallyGeodesicSphereBundle}, there exists a complete totally geodesic immersed submanifold
$\Sigma=(S,i)$ of codimension $k+1$ in $M$ such that if $N\Sigma$ is the $k$-sphere bundle over $\Sigma$ in $M$,
then $\hat{\mathi}$ defines a covering map from $\hat{\Sigma}$ onto $N\Sigma$. The special
Lagrangian angle of $N\Sigma$ equals $k\pi/2=\theta$. Thus, if $\theta\neq(n-1)\pi/2$, we obtain a contradictionand it thus
follows that $\partial\hat{\Cal{F}}_{\rho,\theta}=\emptyset$. On the other hand, if $\theta=(n-1)\pi/2$, then  $k=(n-1)$ and $\Sigma$ is
therefore a geodesic. The result now follows.\qed
\medskip
\noindent Finally we obtain Theorem \procref{ThmRigidity} as a corollary:
\medskip
{\bf\noindent Proof of Theorem \procref{ThmRigidity}:\ }Let $\Sigma=(S,i)$ be a convex, immersed hypersurface in $\Cal{F}_{\rho,\theta}(M)$. Let $(\Sigma_n)_\ninn=(S_n,i_n)_\ninn$
be a sequence of convex immersed submanifolds in $\Cal{F}_{\rho,\theta}(M)$ such that $i_n\rightarrow i$ in the $C^0$ topology. For all $n$, let $\hat{\mathi}_n$ be the
$\rho$-Gauss lifting of $i_n$. The convexity of $\Sigma_n$ for all $n$ and the convergence of $(i_n)_\ninn$ to $i_0$ ensures that no subsequence of $\hat{\mathi}_n$ can converge
to a sphere bundle over a geodesic. It thus follows by Theorem \procref{ThmCompactnessOfHypersurfaces} that every subsequence of $(i_n)_\ninn$ subconverges in the $C^\infty$
Cheeger/Gromov topology (thus modulo reparametrisation) to an immersion $i_0'$. Trivially, $i_0'=i_0$. It thus follows that $(i_n)_\ninn$ converges itself in $C^\infty$ Cheeger/Gromov
topology to $i_0$. The result now follows by Lemma \procref{LemmaLocalRigidity}.\qed
\medskip
\noindent Finally, as remarked in the introduction, the techniques used to prove Theorems \procref{ThmCompactness} and \procref{ThmCompactnessOfHypersurfaces} may be trivially adapted to
yield the following analogous result for Gaussian curvature:
\proclaim{Proposition \nextprocno}
\noindent Let $M$ be a Riemannian manifold of dimension $(n+1)$. Let $(\Sigma_n,p_n)_\ninn$ be a sequence of complete, immersed submanifolds of 
constant Gaussian curvature. For all $n$, let $A_n$ be the shape operator of $\Sigma_n$. If there exists $K>0$ such that, for all $n$:
$$
\|A_n\| \leqslant K,
$$
\noindent then $(\Sigma_n,p_n)_\ninn$ subconverges in the $C^\infty$ Cheeger/Gromov sense to a complete, immersed submanifold of constant Gaussian curvature.
\endproclaim
\proclabel{CompactnessForGaussianCurvature}
\inappendicestrue
\global\headno=0
\newhead{The Maximum Principle}
\noindent In this appendix, we briefly prove a version of the maximum principle required in the proof of Theorem \procref{ThmCompactnessOfHypersurfaces}.
\medskip
\noindent Let $\Omega$ be an open subset of $\Bbb{R}^k$. Let $(g_n)_\ninn$ be a sequence of smooth metrics over 
$\Omega$ converging smoothly over $\Omega$ to $g_0$. For all $n\in\Bbb{N}\munion\left\{0\right\}$, let $\Delta_n$
be the Laplacian of $g_n$. Thus, if $\Gamma_n$ is the Christophel symbol of the Levi-Civita covariant derivative of
$g_n$, then, for every $f:\Omega\rightarrow\Bbb{R}$:
$$
\Delta_n f = (g_n)^{ij}\partial_i\partial_j f - (g_n)^{ij}(\Gamma_n)^k_{ij}\partial_k f.
$$
\noindent Let $(b_n)_\ninn$ be a sequence of smooth vector fields over $\Omega$ converging smoothly over $\Omega$ to $b_0$.
\proclaim{Lemma \nextprocno}
\noindent Let $c$ and $\alpha$ be strictly positive real numbers. For all $n\in\Bbb{N}$, let $f_n:\Omega\rightarrow]0,\infty[$ be a smooth, positive valued function such that:
$$
\Delta_n f_n^\alpha + \langle b_n, df_n^\alpha\rangle - c f_n^\alpha \leqslant 0.
$$
\noindent Suppose that $(f_n)_\ninn$ converges $C^\infty$ to $f_0:\Omega\rightarrow[0,\infty[$.
\medskip
\noindent If there exists $p\in\Omega$ such that $f_0(p)=0$, then $f_0$ is identically zero.
\endproclaim
\proclabel{MaximumPrinciple}
\remark Care is required since the function $f_0^\alpha$ may not be $C^1$ at $0$.
\medskip
\proof Let $B$ be a closed ball about $p$ in $\Omega$. For all $n$, we define the operator $D_n$ by:
$$
D_n f = \Delta_n f + \langle b_n, df\rangle - c f.
$$
\noindent We consider this operator acting on the space of smooth functions over $B$ which vanish on the boundary. Since $c<0$, by the maximum principle (see
\cite{GilbTrud}), this operator has trivial kernel on this space. Thus, by classical Fredholm theory, for any continuous function $\varphi$
on $\partial B$, there exists a unique solution to the Dirichlet problem given by the operator $D_n$ with boundary values equal to $\varphi$. 
In particular, for all $n\in\Bbb{B}\munion\left\{0\right\}$, since $c<0$, there exists a unique function $u_n:B\rightarrow\Bbb{R}$ such that 
$\Delta_n u_n = 0$ and:
$$
u_n|_{\partial B} = f_n^\alpha|_{\partial B}.
$$
\noindent Since $\Delta_n$ converges $C^\infty$ to $\Delta_0$ and $f_n^\alpha$ converges $C^0$ to $f_0^\alpha$, $(u_n)_\ninn$ converges $C^\infty$ 
to $u_0$. For all $n\neq 0$:
$$
\Delta_n (f_n^\alpha - u_n) \leqslant 0.
$$
\noindent Since $f_n^\alpha - u_n$ is smooth, it follows by the maximum principle that $(f_n^\alpha - u_n)$ cannot have a non-positive minimum. Thus, for all $n\neq 0$:
$$
f_n^\alpha \geqslant u_n.
$$
\noindent Taking limits:
$$
f_0^\alpha \geqslant u_0.
$$
\noindent Thus $u_0(p)=0$. Since $u_0\geqslant 0$ along $\partial B$, it follows that $u_0$ has a non-positive minimum in the interior of $B$. Since 
$u_0$ is smooth, the strong maximum principle implies that $u_0$ vanishes identically. Since $u_0$ coincides with $f_0^\alpha$ along 
$\partial B$, and since $B$ is arbitrary, it follows that $f_0$ vanishes identically in a neighbourhood of $p$. The result now follows by a 
standard open/closed argument.\qed
\newhead{Bibliography}
{\leftskip = 5ex \parindent = -5ex
\leavevmode\hbox to 4ex{\hfil\cite{Aronsz}}\hskip 1ex{Aronszajn N., A unique continuation theorem for elliptic differential equations or inequalities of the second order, {\sl J. Math. Pures Appl.} {\bf 36} (1957),
2359--239}
\medskip
\leavevmode\hbox to 4ex{\hfil\cite{BallGromSch}}\hskip 1ex{Ballman W., Gromov M., Schroeder V., {\it Manifolds of nonpositive curvature}, Progress in Mathematics, {\bf 61}, Birkh\"auser, Boston, (1985)}
\medskip
\leavevmode\hbox to 4ex{\hfil\cite{CalabiB}}\hskip 1ex{Calabi E., Improper affine hyperspheres of convex type and a generalisation of a theorem by K. J\"orgens, {\sl Michigan Math. J.} {\bf 5} (1958), 105--126} 
\medskip
\leavevmode\hbox to 4ex{\hfil\cite{Corlette}}\hskip 1ex{Corlette K., Immersions with bounded curvature, {\sl Geom. Dedicata} {\bf 33} (1990), no. 2, 153--161}
\medskip
\leavevmode\hbox to 4ex{\hfil\cite{GilbTrud}}\hskip 1ex{Gilbard D., Trudinger N. S., {\sl Elliptic partical differential equations of second order}, Die Grund\-lehren der mathemathischen Wissenschagten, {\bf 224},
Springer-Verlag, Berlin, New York (1977)}
\medskip
\leavevmode\hbox to 4ex{\hfil\cite{GromB}}\hskip 1ex{Gromov M., Pseudoholomorphic curves in symplectic manifolds, {\sl Invent. Math.} {\bf 82} (1985),  no. 2, 307--347}
\medskip
\leavevmode\hbox to 4ex{\hfil\cite{GromA}}\hskip 1ex{Gromov M., {\sl Metric Structures for Riemannian and Non-Riemannian Spaces}, Progress in Mathematics, {\bf 152}, Birkh\"auser, Boston, (1998)}
\medskip
\leavevmode\hbox to 4ex{\hfil\cite{HarveyLawson}}\hskip 1ex{Harvey R., Lawson H. B. Jr., Calibrated geometries, {\sl Acta. Math.} {\bf 148} (1982), 47--157}
\medskip
\leavevmode\hbox to 4ex{\hfil\cite{Hitchin}}\hskip 1ex{Hitchin N., The moduli space of special Lagrangian submanifolds, {\sl Ann. Scuola Norm. Sup. Pisa Cl. Sci. (4)} {\bf 25} (1997), no. 3-4, 503--515}
\medskip
\leavevmode\hbox to 4ex{\hfil\cite{JostXin}}\hskip 1ex{Jost J., Xin Y.L., A Bernstein theorem for special Lagrangian graphs, {\sl Calc. Var.} {\bf 15} (2002), 299--312}
\medskip
\leavevmode\hbox to 4ex{\hfil\cite{Joyce}}\hskip 1ex{Joyce D., Lectures on special Lagrangian geometry, {\sl Global theory of minimal surfaces}, Clay Math. Proc., Amer. Math. Soc., Providence, RI, (2005), {\bf 2}, 667--695.}
\medskip
\leavevmode\hbox to 4ex{\hfil\cite{LabC}}\hskip 1ex{Labourie F., Probl\`eme de Minkowksi et surfaces \`a courbure constante dans les vari\'e\t'es hyperboliques, {\sl Bull. Soc. Math. France} {\bf 119} (1991), no. 3, 307--325}
\medskip
\leavevmode\hbox to 4ex{\hfil\cite{LabB}}\hskip 1ex{Labourie F., Probl\`emes de Monge-Amp\`ere, courbes holomorphes et laminations, {\sl GAFA} {\bf 7}, no. 3, (1997), 496--534}
\medskip
\leavevmode\hbox to 4ex{\hfil\cite{LabA}}\hskip 1ex{Labourie F., Un lemme de Morse pour les surfaces convexes, {\sl Invent. Math.} {\bf 141} (2000), 239--297}
\medskip
\leavevmode\hbox to 4ex{\hfil\cite{Peterson}}\hskip 1ex{Peterson P., {\sl Riemannian Geometry}, Graduate Texts in Mathematics, {\bf 171}, Springer Verlag, New York, (1998)}
\medskip
\leavevmode\hbox to 4ex{\hfil\cite{Pog}}\hskip 1ex{Pogorelov A. V., On the improper convex affine hyperspheres, {\sl Geometriae Dedicata} {\bf 1} (1972), 33--46}
\medskip
\leavevmode\hbox to 4ex{\hfil\cite{RosSpruck}}\hskip 1ex{Rosenberg H., Spruck J. On the existence of convex hyperspheres of constant Gauss curvature in hyperbolic space, {\sl J. Diff. Geom.} {\bf 40} (1994), no. 2,
379--409} 
\medskip
\leavevmode\hbox to 4ex{\hfil\cite{SalMcDuffB}}\hskip 1ex{McDuff D., Salamon D., {\sl $J$-holomorphic curves and quantum cohomology}, University Lecture Series, {\bf 6}, AMS, Providence, (1994)}
\medskip
\leavevmode\hbox to 4ex{\hfil\cite{SalMcDuffA}}\hskip 1ex{McDuff D., Salamon D., {\sl Introduction to symplectic topology}, Oxford, (1995)}
\medskip
\leavevmode\hbox to 4ex{\hfil\cite{SmiE}}\hskip 1ex{Smith G., Th\`ese de doctorat, Paris (2004)}%
\medskip
\leavevmode\hbox to 4ex{\hfil\cite{SmiG}}\hskip 1ex{Smith G., Pointed $k$-surfaces, {\sl Bull. Soc. Math. France} {\bf 134} (2006), no. 4, 509--557}
\medskip
\leavevmode\hbox to 4ex{\hfil\cite{SmiF}}\hskip 1ex{Smith G., An Arzela-Ascoli theorem for immersed submanifolds, to appear in {\sl Ann. Fac. Sci. Toulouse Math.\/}, math.DG/0510231}\par%
\medskip
\leavevmode\hbox to 4ex{\hfil\cite{SmiH}}\hskip 1ex{Smith G., Constant special Lagrangian hypersurfaces in quasi-Fuchsian hyperbolic ends, in preparation}
\medskip
\leavevmode\hbox to 4ex{\hfil\cite{Smoczyk}}\hskip 1ex{Smoczyk K., Longtime existence of the Lagrangian mean curvature flow, {\sl Calc. Var. Partial Differential Equations} {\bf 20} (2004), no. 1, 25--46}
\medskip
\leavevmode\hbox to 4ex{\hfil\cite{Yuan}}\hskip 1ex{Yuan Y., A Bernstein problem for special Lagrangian equations, {\sl Invent. Math.} {\bf 150} (2002), 117--125}\par%
}%
\enddocument

%% file: references.tex
\global\def\_@citation@Aronsz{1}
\global\def\_@citation@BallGromSch{2}
\global\def\_@citation@CalabiB{3}
\global\def\_@citation@Corlette{4}
\global\def\_@citation@GilbTrud{5}
\global\def\_@citation@GromB{6}
\global\def\_@citation@GromA{7}
\global\def\_@citation@HarveyLawson{8}
\global\def\_@citation@Hitchin{9}
\global\def\_@citation@JostXin{10}
\global\def\_@citation@Joyce{11}
\global\def\_@citation@LabC{12}
\global\def\_@citation@LabB{13}
\global\def\_@citation@LabA{14}
\global\def\_@citation@Peterson{15}
\global\def\_@citation@Pog{16}
\global\def\_@citation@RosSpruck{17}
\global\def\_@citation@SalMcDuffB{18}
\global\def\_@citation@SalMcDuffA{19}
\global\def\_@citation@SmiE{20}
\global\def\_@citation@SmiG{21}
\global\def\_@citation@SmiF{22}
\global\def\_@citation@SmiH{23}
\global\def\_@citation@Smoczyk{24}
\global\def\_@citation@Yuan{25}
\global\def\_@proc@ThmCompactness{1.2}
\global\def\_@proc@ThmRigidity{1.3}
\global\def\_@proc@ThmCompactnessOfHypersurfaces{1.4}
\global\def\_@subhead@ImmersedSubmanifoldsAndTheCheegerGromovTopology{2.1}
\global\def\_@head@HeadStructuresSLClassiques{2}
\global\def\_@subhead@HeadStructuresSLClassiques{2.2}
\global\def\_@proc@LemmaGaussLiftingLegendrian{2.1}
\global\def\_@proc@EquivalenceSLCurvatureAndSLImmersions{2.2}
\global\def\_@proc@ThmUniqueContinuation{2.3}
\global\def\_@proc@LemmaInfinitesimalDeformation{3.1}
\global\def\_@proc@LemmeBijectivite{3.2}
\global\def\_@proc@LemmaLocalRigidity{3.3}
\global\def\_@head@HeadCompacite{4}
\global\def\_@proc@ThmYuanPlanes{4.1}
\global\def\_@proc@LemmeArzelaAscoliCorlette{4.3}
\global\def\_@proc@LemmeDuLambdaMaximum{4.4}
\global\def\_@proc@LemmaCompactness{4.5}
\global\def\_@proc@LemmeRegulariteElliptique{4.6}
\global\def\_@proc@DALemmaSymmOfDA{5.1}
\global\def\_@proc@DALemmaSymmOfDDA{5.1}
\global\def\_@proc@ConstSLLemmaDA{5.2}
\global\def\_@subhead@CovariantDerivatives{5.2}
\global\def\_@proc@CDLapLemmaChristSymbol{5.3}
\global\def\_@proc@CDLapLemmaLaplacian{5.4}
\global\def\_@subhead@UsefulFunction{5.3}
\global\def\_@proc@UsefulLemmaSymmOfDDA{5.5}
\global\def\_@proc@DALemmaRelationAandTildeA{5.6}
\global\def\_@proc@LemmaLapOfLog{5.7}
\global\def\_@proc@UsefulLemmaLaplaceF{5.8}
\global\def\_@proc@Subharmonic{5.9}
\global\def\_@proc@LemmaBoundaryVerticalForThetaNotInteger{5.11}
\global\def\_@proc@GCSLemmaDecomposition{5.12}
\global\def\_@proc@LemmaLimitIsTotallyGeodesicSphereBundle{5.13}
\global\def\_@proc@CompactnessForGaussianCurvature{5.14}
\global\def\_@proc@MaximumPrinciple{A.1}